\documentclass[preprint,12pt]{elsarticle}
\usepackage{amssymb}
\usepackage{amsmath}
\usepackage{amsfonts}
\begin{document}

\begin{frontmatter}

\title{Cone-theoretic generalization of total positivity}
\author{O.Y. Kushel}
\ead{kushel@mail.ru}
\address{Institut f\"{u}r Mathematik, MA 4-5, Technische Universit\"{a}t Berlin, D-10623 Berlin, Germany}

\begin{abstract}
This paper is devoted to the generalization of the theory of total positivity.
We say that a linear operator $A: {\mathbb R}^n \rightarrow {\mathbb R}^n$ is generalized totally positive (GTP), if its $j$th exterior power $\wedge^j A$ preserves a proper cone $K_j \subset \wedge^j {\mathbb R}^n$ for every $j = 1, \ \ldots, \ n$. We also define generalized strictly totally positive (GSTP) operators. We prove that the spectrum of a GSTP operator is positive and simple, moreover, its eigenvectors are localized in special sets. The existence of invariant cones of finite ranks is shown under some additional conditions. Some new insights and alternative proofs of the well-known results of Gantmacher and Krein describing the properties of TP and STP matrices are presented.

\end{abstract}

\begin{keyword}
Cones of rank $k$ \sep Cone-preserving maps \sep Gantmacher--Krein theorem \sep Total
positivity \sep Compound matrices \sep Exterior products.

 \MSC Primary 15A48 \sep Secondary 15A18 \sep 15A75
\end{keyword}

\end{frontmatter}

\newtheorem{thm}{Theorem}
\newtheorem{lem}[thm]{Lemma}
\newtheorem{prop}[thm]{Proposition}
\newtheorem{cor}[thm]{Corollary}
\newdefinition{rmk}{Remark}
\newproof{pf}{Proof}
\newproof{pot}{Proof of Theorem \ref{thm2}}


\section{Introduction}

The theory of totally positive matrices and kernels
started with Kellog \cite{KEL1,KEL2} and mainly developed in monographs \cite{GANT} by Gantmacher and Krein and \cite{KARL1} by Karlin, nowadays becomes an interesting and important part of the modern analysis. A matrix ${\mathbf A}$ is called
{\it positive} if all its elements $a_{ij}$
are positive. A $n \times n$ matrix ${\mathbf A}$ is called {\it strictly totally positive (STP)} if its $j$th compound matrix ${\mathbf A}^{(j)}$ is positive for every $j = 1, \ \ldots, \ n$.
(Recall that ${\mathbf A}^{(j)}$ is the matrix that consists of all the minors
$A\begin{pmatrix}
  i_1 & \ldots & i_j \\
  k_1 & \ldots & k_j
\end{pmatrix}$, where $1 \leq i_1 < \ldots < i_j \leq n, \ 1 \leq \ k_1 < \ldots < k_j \leq n$,  of the initial matrix ${\mathbf A}$. The minors are listed in the lexicographic order. The matrix ${\mathbf A}^{(j)}$ is $\binom{n}j \times \binom{n}j$ dimensional, where $\binom{n}j = \dfrac{n!}{j!(n-j)!}$. The first compound matrix ${\mathbf A}^{(1)}$ is equal to ${\mathbf A}$).

We introduce the following definition which gives a natural generalization of the class of STP matrices. Given a family of proper cones $\{K_1, \ \ldots, \ K_n\}$, $K_j \subset {\mathbb R}^{\binom{n}j}$, we call a $n \times n$ matrix ${\mathbf A}$ {\it generalized strictly totally positive (GSTP)} with respect to $\{K_1, \ \ldots, \ K_n\}$ if its $j$th compound matrix ${\mathbf A}^{(j)}$ maps $K_j \setminus \{0\}$ into ${\rm int}(K_j)$ for every $j = 1, \ \ldots, \ n$.

The following result of Schoenberg is known for STP matrices (see, for example, \cite{PINK}). Let us recall the following two ways of counting for the number of sign changes of a vector $x = (x^1, \ \ldots, \ x^n) \in {\mathbb R}^n$. $S^-(x)$ denotes the number of sign changes in the sequence $(x^1, \ \ldots, \ x^n)$ of the coordinates with zero terms discarded. $S^+(x)$ denotes the maximum number of sign changes in the sequence $(x^1, \ \ldots, \ x^n)$ where zero terms are arbitrarily assigned values $\pm 1$ (see, for example, \cite{PINK}, p. 76).

\begin{thm}[Schoenberg]
 Let a $n \times n$ matrix ${\mathbf A}$ be STP. Then the following inequality holds for each nonzero vector $x \in {\mathbb R}^n$:
$$S^+({\mathbf A}x) \leq S^-(x).$$
\end{thm}

We construct special sets $T(K_j) \subseteq {\mathbb R}^n$ with respect to the cones $K_j \subset {\mathbb R}^{\binom{n}j}$. Thus we obtain the following generalization of the Schoenberg theorem.
\setcounter{thm}{24}
\begin{thm}
Let a $n \times n$ matrix ${\mathbf A}$ be GSTP with respect to $\{K_1, \ \ldots, \ K_n\}$. Then the interior of the set $T(K_j)$ is nonempty and the inclusion $x \in T(K_j) \setminus \{0\}$ implies the inclusion ${\mathbf A}x \in {\rm int}(T(K_j))$ for every $j = 1, \ldots, \ n$.
\end{thm}

We also generalize the classical Gantmacher-Krein theorem (see, for example, \cite{AN,PINK,PIN1}) to the case of GSTP matrices.
\setcounter{thm}{1}
\begin{thm}[Gantmacher, Krein]
Let the matrix ${\mathbf A}$ of a linear operator $A: {\mathbb R}^n
\rightarrow {\mathbb R}^n$ be STP. Then all the
eigenvalues of the operator $A$ are positive and simple:
$$\rho(A)= \lambda_1 > \lambda_2 > \ldots > \lambda_n > 0.$$
 The first eigenvector corresponding to the maximal
eigenvalue $\lambda_1$ is strictly positive and the $j$th
eigenvector $x_j$ corresponding to the $j$th in absolute value
eigenvalue $\lambda_j$ has exactly $j - 1$ changes of sign. Moreover, the following inequalities hold:
$$q-1 \leq S^-(\sum_{i=q}^p c_i x_i) \leq S^+(\sum_{i=q}^p c_i x_i) \leq p-1$$
for each $1 \leq q \leq p \leq n$ and $\sum\limits_{i=q}^p c_i^2 \neq 0$.
\end{thm}

We construct special sets $T(K_1, \ldots,
K_j) \subseteq {\mathbb R}^n$, $j = 1, \ \ldots, \ n$ with respect to the given family of cones $\{K_1, \ \ldots, \ K_n\}$, $K_j \subset {\mathbb R}^{\binom{n}j}$.
\setcounter{thm}{21}
\begin{thm} Let a linear operator $A: {\mathbb{R}}^{n}
\rightarrow {\mathbb{R}}^{n}$ be GSTP with respect to a totally positive structure $\{K_1, \ \ldots, \ K_n\}$. Then all the eigenvalues of the operator $A$ are
positive and simple:
$$\rho(A)= \lambda_1 > \lambda_2 > \ldots > \lambda_n > 0.$$
 The first eigenvector $x_1$ corresponding to the maximal
eigenvalue $\lambda_1$ belongs to ${\rm int} (K_1)$ and the
$j$th eigenvector $x_j$ corresponding to the $j$th in absolute value
eigenvalue $\lambda_j$ belongs to ${\rm int}(T(K_1, \ldots,
K_j))\setminus T(K_1, \ldots, K_{j-1})$. Moreover, the following
inclusions hold:
$$\sum_{i=q}^p c_i x_i \in {\rm int}(T(K_1, \ \ldots, \
K_p))\setminus T(K_1, \ \ldots, \ K_{q - 1})$$ for each $1 \leq q \leq p
\leq n$ and $c_p \neq 0$;
$$\sum_{i=q}^p c_i x_i \in \overline{T(K_1, \ \ldots, \
K_p)\setminus T(K_1, \ \ldots, \ K_{q - 1})}$$ for each $1 \leq q \leq p
\leq n$.
\end{thm}
\setcounter{thm}{2}
The organization of this paper is as follows. In Section 2, we introduce basic definitions concerning exterior powers of finite-dimensional spaces. In Section 3, we recall basic definitions of the theory of cones and provide some examples. Here we also give the definition of a cone of a finite rank. Section 4 deals with a certain duality between cones of rank $j$ in ${\mathbb R}^n$ and proper cones in its $j$th exterior power $\wedge^j {\mathbb R}^n$. In particular, we construct a special set $T(K_j) \subset {\mathbb R}^n$ for a given proper cone $K_j \subset \wedge^j {\mathbb R}^n$ and study its topological properties. We conclude that under certain additional assumptions the set $T(K_j)$ is a cone of rank $j$. In Section 5, we construct a set $T(K_1,\ldots,K_j)$ with respect to a family of proper cones $K_1 \subset {\mathbb R}^n$, $K_2 \subset \wedge^2 {\mathbb R}^n$, $\ldots$, $K_j \subset \wedge^j {\mathbb R}^n$. We give examples of such sets and study their topological properties. In Section 6, we state the main results of the theory of cone-preserving maps. In Section 7, we recall basic facts concerning exterior powers of linear operators in ${\mathbb R}^n$. In Section 8, we introduce the concepts of generalized total positivity (with respect to a given family of proper cones), generalized strict total positivity and generalized sign-regularity. Such definitions provides natural generalizations of the classes of totally positive, strictly totally positive and sign-regular matrices, respectively. Basic properties of GTP, GSTP and GSR operators are listed in Section 9. The results of this section shows that the class of GSR operators covers the entire class of operators with real spectrum. In Section 10, we state and prove the generalization of the result of Schoenberg concerning variation-diminishing properties of SR matrices. The results of this section shows that a GSR (with respect to a family of proper cones) operator preserves conic sets constructed as it was shown in Sections 4-5. Our main result concerning spectral properties of GSTP operators is proved in Section 11. In this section we also provide some conditions for a family of cones which are necessary for the existence of at least one GSTP operator. Then we state and prove a stronger statement describing invariant sets of a GSSR operator. In Section 12, we deduce the classical results on TP and STP operators (which are special cases of GTP and GSTP operators) from the preceding reasoning. In Section 13, we study one more special case of GTP matrices, in particular, matrices every compound of which is diagonally similar to a positive matrix. We list some special properties of such matrices and provide examples which shows that this statements are not valid for arbitrary GTP operators. Some conclusions are given in Section 14.

\section{Exterior powers of the space ${\mathbb R}^n$}
Let ${\mathbb R}^n$ denote $n$-dimensional Euclidean space, and
$({\mathbb R}^n)'$ denote an adjoint space of all linear
functionals on ${\mathbb R}^n$. Since $({\mathbb R}^n)'$ is also
$n$-dimensional, we consider linear functionals from $({\mathbb R}^n)'$
as vectors from ${\mathbb R}^n$.

 Let us recall some basic
definitions and statements about the tensor and exterior powers
of the space ${\mathbb R}^n$ (for more complete information see
\cite{GLALU, TSOY, YOK}).

Let $j = 2, \ \ldots, \ n$. The space of all multilinear functionals on $
\times^j ({\mathbb R}^n)^ \prime$ is called the {\it $j$th tensor power} of
the space ${\mathbb R}^n$ and denoted by $ \otimes^j {\mathbb
R}^n$. Its elements are called {\it tensors}.

 Let $x_1, \ \ldots, \ x_j$ be arbitrary vectors from ${\mathbb R}^n$. Then the multilinear
functional $x_1 \otimes \ldots \otimes x_j: (\times^j ({\mathbb
R}^n)^ \prime) \rightarrow {\mathbb R}$ which acts according to
the rule
$$(x_1 \otimes \ldots \otimes x_j)(f_1, \ldots, f_j)=
\langle x_1,f_1 \rangle \ldots \langle x_j,f_j \rangle,$$ is called a {\it tensor
product} of the vectors $x_1, \ \ldots, \ x_j$. (Here the linear functionals
 $f_1, \ \ldots, \ f_j \in ({\mathbb R}^n)^ \prime$ are considered as
vectors from ${\mathbb R}^n$).

 The $j$-th tensor power $ \otimes^j {\mathbb R}^n$ of the space ${\mathbb R}^n$
is spanned by elementary tensor
products of the form $x_1 \otimes \ldots \otimes x_j$ where $x_1, \ \ldots, \ x_j  \in {\mathbb
R}^n$. Examine an arbitrary basis
$e_1, \ \ldots, \ e_n$ in ${\mathbb R}^n$.  Then all the possible tensor products
of the form $e_{i_1} \otimes \ldots \otimes e_{i_j}$ $ \ (1 \leq i_1, \ \ldots, \ i_j \leq n)$ of
the initial basic vectors form a basis in $ \otimes^j {\mathbb
R}^n$. It follows that the space $\otimes^j {\mathbb R}^n$ is finite-dimensional with
${\rm dim}(\otimes^j {\mathbb R}^n) = n^j$.

Let $(i_1, \ \ldots, \ i_j)$ be a permutation of the set $[j] = \{1, \
\ldots, \ j\}$. Define
$$\chi(i_1, \ldots, i_j) = \left\{\begin{array}{cc} 1,& \mbox{if
the permutation $(i_1, \ \ldots, \ i_j)$ is even};
\\[10pt] - 1, & \mbox{if the permutation is odd.}\end{array}\right.$$

{\it The $j$th exterior power} $ \wedge^j {\mathbb
R}^n$ of the space ${\mathbb R}^n$ is a subspace of the space $ \otimes^j {\mathbb
R}^n$ consisting of all antisymmetric tensors (i.e. all the tensors
$\varphi$ for which $\varphi(f_1, \ldots, f_j) = \chi(i_1, \ldots, i_j)
\varphi(f_{i_1}, \ldots, f_{i_j})$ where
$f_1, \ldots, f_j$ are arbitrary functionals from $({\mathbb R}^n)^\prime$).

 Let $x_1, \ldots, \ x_j$ be arbitrary
vectors from ${\mathbb R}^n$. Then the multilinear functional $x_1
\wedge \ldots \wedge x_j: \times^j ({\mathbb R}^n)^ \prime
\rightarrow {\mathbb R}$ which acts according to the rule
$$(x_1 \wedge \ldots \wedge x_j)(f_1, \ldots, f_j)=
\sum_{(i_{1}, \ldots, i_j)} \chi(i_1, \ldots, i_j)
(x_{i_1} \otimes \ldots \otimes x_{i_j})(f_1, \ldots, f_j) = $$ $$ = \sum_{(i_{1}, \ldots, i_j)} \chi(i_1, \ldots, i_j) \ \langle
x_{i_1}, f_1 \rangle \ldots \langle x_{i_j}, f_j \rangle$$ is
called an {\it exterior product} of the vectors $x_1, \ \ldots, \
x_j$. Here the sum is taken with respect to all the permutations
$(i_{1}, \ \ldots, \ i_j)$ of $[j]$ and linear functionals  $f_1, \ \ldots, \ f_j \in ({\mathbb
R}^n)^ \prime$ are considered as vectors from ${\mathbb R}^n$.

It is easy to see, that the exterior product $x_1 \wedge \ldots \wedge x_j$ is
antisymmetric, i.e. the following equality holds for every permutation $(i_1, \ \ldots, \
i_j)$ of $[j]$:
$$x_{i_1}\wedge\ldots \wedge x_{i_j} = \chi(i_1, \ \ldots, i_j)
(x_1 \wedge \ldots \wedge x_j).$$

The space $\wedge^j {\mathbb R}^n$ is spanned by all the exterior products $x_1 \wedge \ldots \wedge x_j$ where $x_1, \
\ldots, \ x_j \in {\mathbb R}^n$. If the vectors $e_1, \ \ldots, \ e_n$ form a basis in the initial space ${\mathbb R}^n$ then the set of all exterior
products of the type $\{e_{i_1} \wedge \ldots \wedge e_{i_j}\}$ where
$1 \leq i_1 < \ldots < i_j \leq n$ forms a canonical basis in the
space $\wedge^j {\mathbb R}^n$ (see \cite{GANT,PINK,POST2}). Thus the space $\wedge^j {\mathbb R}^n$ is finite dimensional with ${\rm dim}(\wedge^j {\mathbb R}^n) = \binom{n}j$. (Here $\binom{n}j = \frac{n!}{j!(n-j)!}$).

A scalar product on $\wedge^j {\mathbb R}^n$ is defined by the formula:
$$\langle x_1 \wedge \ldots \wedge x_j, \ y_1 \wedge \ldots \wedge y_j\rangle = (x_1 \wedge \ldots \wedge x_j)(y_1, \ldots, y_j) = $$ $$ =\sum_{(i_{1}, \ldots, i_j)} \chi(i_1, \ldots, i_j) \ \langle
x_{i_1}, y_1 \rangle \ldots \langle x_{i_j}, y_j \rangle.$$
 It follows that the adjoint space $(\wedge^j{\mathbb R}^n)'$ can be considered as $\wedge^j({\mathbb R}^n)'$ (see \cite{YOK}, p. 88).

Let the element $\varphi \in \wedge^j {\mathbb R}^n$ be represented in the form of the exterior product $x_1 \wedge \ldots \wedge x_j$ of some vectors $x_1, \ \ldots, \ x_j \in {\mathbb R}^n$.  Then $\varphi$ is called a {\it simple $j$-vector}. The set of all simple $j$-vectors is called the {\it Grassmann cone} and denoted $\barwedge^j {\mathbb R}^n$. The equality $\barwedge^j {\mathbb R}^n = \wedge^j {\mathbb R}^n$ holds only for $j = 1, n-1$ and $n$. (Note that $\wedge^1 {\mathbb R}^n = {\mathbb R}^n$ and $\wedge^n{\mathbb R}^n = {\mathbb R}$). If $j = 2, \ \ldots, \ n-2$, then we can find elements of $\wedge^j {\mathbb R}^n$ which can not be represented as simple $j$-vectors (see \cite{POST2}, p. 83). It is not difficult to see that {\it the set $\barwedge^j {\mathbb R}^n$ is uniform (i.e. the equality $\alpha \barwedge^j {\mathbb R}^n = \barwedge^j {\mathbb R}^n$ is true for every nonzero $\alpha \in {\mathbb R}$) and closed in the space $\wedge^j {\mathbb R}^n$.}

Let us define a map ${\mathcal A}_j$ acting from the set of
all $j$-dimensional subspaces of ${\mathbb R}^n$ to the set of
$1$-dimensional subspaces (i.e. lines) of $\wedge^j {\mathbb R}^n$
according to the following rule:
$$ {\mathcal A}_j(L) = \{t(x_1 \wedge \ldots \wedge x_j)\}_{t \in {\mathbb R}}, $$
where $L$ is a $j$-dimensional subspace of ${\mathbb R}^n$,
$x_1, \ \ldots, \ x_j$ are $j$ arbitrary linearly independent
vectors from $L$.

It is not difficult to see that the map ${\mathcal A}_j$ is well-defined,
i.e. {\it if $x_1, \ \ldots, \ x_j$ and $y_1, \ \ldots, \ y_j$ are
two sets of linearly independent vectors, which belongs to the
same $j$-dimensional subspace $L$, then their exterior products
$x_1 \wedge \ldots \wedge x_j$ and $y_1 \wedge \ldots \wedge y_j$
are collinear} (see, e.g., \cite{POST2}).

The map ${\mathcal A}_j$ is a bijective map between all $j$-dimensional subspaces of ${\mathbb R}^n$ and all lines of $\barwedge^j {\mathbb R}^n$ (see \cite{POST2}, p. 86).

Let us consider the $(n-1)th$ exterior power of the $n$-dimensional space ${\mathbb R}^n$. Note that $\binom{n}{n-1} = n$, thus ${\rm dim}(\wedge^{n-1} {\mathbb R}^n) = n$. All the exterior products of the type $\{e_{i_1}\wedge \ldots \wedge e_{i_{n-1}}\}$ where $1 \leq i_1
< \ldots < i_{n-1} \leq n$ of the initial basic vectors $e_1, \
\ldots, \ e_n$ form a basis in $ \wedge^{n-1} {\mathbb R}^n$. Let us define a
bijective linear operator ${\mathcal J}_n: \wedge^{n-1} {\mathbb
R}^n \rightarrow {\mathbb R}^n$ in the following way:
$${\mathcal J}_n (e_{i_1}\wedge \ldots \wedge e_{i_{n-1}}) = (-1)^{k+1}e_k,$$
where $k = [n]\setminus\{i_1, \ \ldots, \
i_{n-1}\}$.

Let $x_i = (x_i^1, \ \ldots, \ x_i^n)$ where $i = 1, \ \ldots, \
n-1$ be $n-1$ arbitrary linearly independent vectors from ${\mathbb R}^n$. Write the
exterior product $x_1 \wedge \ldots \wedge x_{n-1}$ in the form $$
x_1 \wedge \ldots \wedge x_{n-1} =\sum_{(i_1, \ldots, i_{n-1})}
\begin{vmatrix}
  x_1^{i_1} & \ldots & x_1^{i_{n-1}} \\
  \ldots & \ldots & \ldots \\
 x_{n-1}^{i_1} & \ldots & x_{n-1}^{i_{n-1}} \\
\end{vmatrix} (e_{i_1}\wedge \ldots \wedge e_{i_{n-1}}).$$

It is not difficult to see, that the vector $${\mathcal J}(x_1 \wedge \ldots
\wedge x_{n-1}) = \sum_{k = 1}^n
\begin{vmatrix}
  x_1^{i_1} & \ldots & x_1^{i_{n-1}} \\
  \ldots & \ldots & \ldots \\
 x_{n-1}^{i_1} & \ldots & x_{n-1}^{i_{n-1}} \\
\end{vmatrix} (-1)^{k+1} e_k =
\begin{vmatrix}
 e_1 & \ldots & e_n \\
  x_{1}^1 & \ldots & x_{1}^n \\
  \ldots & \ldots & \ldots \\
 x_{n-1}^1 & \ldots & x_{n-1}^n \\
\end{vmatrix}$$ is
orthogonal to the hyperplane spanned by the vectors $x_1, \
\ldots, \ x_{n-1}$.

\section{ Conic sets: basic definitions and statements}
Let us recall some basic definitions of the theory of cones (see
\cite{BERPL,KRASN1,SOB2,TAM}).

 A closed subset $K
\subset {\mathbb R}^n$ is called a {\it proper cone}, if it is a
convex cone (i.e. for any $x, y \in K, \ \alpha \geq 0$ we have
$x+y, \ \alpha x \in K$), pointed ($K \cap (-K) = \{0\}$) and solid
(${\rm int}(K) \neq \emptyset$).

The set $K^* \subset ({\mathbb R}^n)'$ defined in the following way
$$ K^* = \{x^* \in ({\mathbb R}^n)': \ \forall y \in K \ \ \langle y, x^* \rangle \geq 0\},$$ is called the {\it adjoint cone} to the cone $K$. The set $K$ is a proper cone in ${\mathbb R}^n$ if and only if $K^*$ is a proper cone in $({\mathbb R}^n)'$. The interior of $K^*$ is defined by the equality $$ {\rm int}(K^*) = \{x^* \in ({\mathbb R}^n)': \ \forall y \in K  \ \ \langle y, x^* \rangle > 0\}.$$

{\bf Example 1.} Let $x_1, \ \ldots, \ x_n \in {\mathbb
R}^n$ be linearly independent vectors. The set $$K = \{\sum_{i =
1}^n c_ix_i : \ c_1, \ \ldots, \ c_n\geq 0\}$$ of all linear
combinations of the vectors $x_1, \ \ldots, \ x_n$ with
nonnegative coefficients is a proper cone in ${\mathbb R}^n$. Such
a cone is called {\it spanned} by the vectors $x_1, \ \ldots, \
x_n$. The cone spanned by the basic vectors $e_1, \ \ldots, \
e_n$ is denoted by ${\mathbb R}^n_+$.

{\bf Example 2.} The set $$K = \{x = (x^1, \ \ldots, \ x^n) \in  {\mathbb R}^n;$$ $$ \sqrt{(x^1)^2 + \ldots + (x^{k-1})^2 + (x^{k+1})^2 + \ldots + (x^n)^2} \leq x^k \}$$ is a proper cone in ${\mathbb R}^n$. Such a cone is called an {\it ice-cream cone}.

{\bf Example 3.} Let $F \subset {\mathbb R}^n$ be a closed, convex and bounded set. The set $K(F)$ of all elements of the form $\alpha x$ where $\alpha \geq 0$, $x \in F$ is a pointed convex cone in ${\mathbb R}^n$. If ${\rm int}(F) \neq \emptyset$, then the cone $K(F)$ is obviously proper.

 Let us examine the cones spanned by the
vectors $\epsilon_1e_1, \ \ldots, \ \epsilon_n e_n$ where each $\epsilon_i$ $(i = 1, \ \ldots, \ n)$ is equal to $+1$ or $-1$. This cone is called a {\it basic cone}.
The space ${\mathbb R}^n$ with a fixed basis $e_1, \ \ldots, \
e_n$ consists of $2^n$ basic cones, one of which is ${\mathbb
R}^n_+$.

We list some properties of basic cones which will be used later.

{\bf 1.} The projection of any basic cone on any basic subspace (i.e on a subspace spanned by any subsystem of the initial basic vectors) is a basic cone in this subspace.

{\bf 2.} If $K$ is a basic cone in
${\mathbb R}^n$, then the adjoint cone $K^*$ is also a basic cone in $({\mathbb R}^n)'$.

As it was mentioned above, every basis $e_1, \ \ldots, \
e_n$ in ${\mathbb R}^n$ defines a basis in the
space $\wedge^j {\mathbb R}^n = \binom{n}j$ which consists of all exterior
products of the form $\{e_{i_1} \wedge \ldots \wedge e_{i_j}\}$, where
$1 \leq i_1 < \ldots < i_j \leq n$. Denote the cone spanned by this exterior basic vectors by $\wedge^j {\mathbb R}^n_+$. Let us call a cone in $\wedge^j {\mathbb R}^n$ spanned by the simple $j$-vectors of the form $\pm (e_{i_1} \wedge \ldots \wedge e_{i_j})$ where $1 \leq i_1 < \ldots < i_j \leq n$ an {\it exterior basic cone} defined by the basis $e_1, \ \ldots, \
e_n$. It is easy to see, that not every basic cone in $\wedge^j {\mathbb R}^n$ is an exterior basic cone.

We list some obvious properties of exterior basic cones.

{\bf 1.} Let $L$ be any basic subspace of ${\mathbb R}^n$. Then the projection of any exterior basic cone on the subspace $\wedge^j L$ of the space $\wedge^j {\mathbb R}^n$ is an exterior basic cone in this subspace.

{\bf 2.} If $K_j$ is an exterior basic cone in
$\wedge^j {\mathbb R}^n$, then the adjoint cone $K_j^*$ is an exterior basic cone in $(\wedge^j{\mathbb R}^n)' = \wedge^j({\mathbb R}^n)'$.

Let us recall the following characterization of a proper cone $K$
(see, for example, \cite{IUS}).

 The angle $\theta_{\max}(K)$ defined by the equality
$$\theta_{\max}(K) = \sup_{x,y \in K\cap S_n}\arccos\langle x,y \rangle, $$
where $S_n$ is the unit sphere in ${\mathbb R}^n$, is called the
{\it maximal angle} of the cone $K$.

Any basic cone in ${\mathbb R}^n$ can be converted using some linear transformation to the cone ${\mathbb
R}^n_+$. Thus we can assume without loss of generality that any basic cone $K$ satisfies the inequality $\theta_{\max}(K) \leq \dfrac{\pi}{2}$.

Besides cones we shall be interested in some other sets in ${\mathbb R}^n$.
Recall the definitions of the following conic sets (see \cite{KRASN1}).

 A closed subset $T \subset {\mathbb
R}^n$ is called {\it a cone of rank $k$} $\ (0 \leq k \leq n)$ if
for every $x \in T$, $\alpha \in {\mathbb R}$ the element $\alpha
x \in T$ and there is at least one $k$-dimensional subspace and no
higher dimensional subspaces in $T$.

For the definition and examples of cones of rank $k$ see also
\cite{KRASN2, SOB1, SOB2}. Note that a cone of rank $k$ is usually not
convex.

{\bf Example 1.} Let $L_1, \ \ldots, \ L_m$ be
subspaces of ${\mathbb R}^n$ with $\max\limits_i \dim (L_i) = k$. Then $\bigcup\limits_{i = 1}^m L_i$ is a
cone of rank $k$ in ${\mathbb R}^n$.

{\bf Example 2.} Let $K \subset {\mathbb R}^n$ be a proper cone.
Then $K \cup (- K)$ is a cone of rank $1$ in ${\mathbb R}^n$, and ${\mathbb R}^n \setminus ({\rm int}(K) \cup {\rm int}(-K))$
is a cone of rank $n-1$ in ${\mathbb R}^n$.

\section{Set $T(K_j)$ and its properties}

Given a proper cone $K_j \subset \wedge^j {\mathbb R}^n = {\mathbb R}^{\binom{n}j}$, $j = 2, \
\ldots, \ n$. Let us define the set $T(K_j) \subset
{\mathbb R}^n$ in the following way:

$$T(K_j) = \overline{\{x_1 \in {\mathbb{R}}^{n}: \exists  \
x_2, \ \ldots, \ x_j \in {\mathbb{R}}^{n}, \ \mbox{for which}}$$ $$\overline{
x_1 \wedge x_2 \wedge \ldots \wedge x_j \in
({\rm int}(K_j) \cup
{\rm int}(- K_j))\}}.$$

Let us define the set $\widehat{T}(K_j)$ in the following way:

$$\widehat{T}(K_j) = \{x_1 \in {\mathbb{R}}^{n}: \exists  \
x_2, \ \ldots, \ x_j \in {\mathbb{R}}^{n}, \ \mbox{for which}$$ $$ \
x_1 \wedge x_2 \wedge \ldots \wedge x_j \in (K_j\cup (- K_j))
\setminus \{0\}\} \cup\{0\}.$$

It is not difficult to see, that the sets $T(K_j)$ and $\widehat{T}(K_j)$ may not coincide for an arbitrary proper cone $K_j \subset \wedge^j {\mathbb R}^n$.

 The following lemma describes the structure of the sets $T(K_j)$ and
$\widehat{T}(K_j)$.

\begin{lem} Let $K_j \subset \wedge^j {\mathbb R}^n$ be
a proper cone. Then the set $\widehat{T}(K_j)$, if it is not $\{0\}$, coincides with the
 set of all $j$-dimensional subspaces $L \subset
{\mathbb R}^n$ for which corresponding lines ${\mathcal A}_j(L)$
belong to $K_j\cup (- K_j)$. The set $T(K_j)$, if it is not empty, coincides with the closure of the set of all $j$-dimensional subspaces $L \subset
{\mathbb R}^n$ for which corresponding lines ${\mathcal A}_j(L)$
belong to ${\rm int}(K_j) \cup
{\rm int}(- K_j)$.\end{lem}

\begin{pf}
$\Leftarrow$ The inclusion $0 \in \widehat{T}(K_j)$ follows from the definition of the set $\widehat{T}(K_j)$. Let an arbitrary nonzero vector $x_1$ belong to a $j$-dimensional subspace $L$ for which the corresponding line ${\mathcal
A}_j(L)$ belongs to $K_j\cup (- K_j)$. Let us show that
$x_1 \in \widehat{T}(K_j)$. Indeed, let us find vectors $
x_2, \ \ldots, \ x_j$ such that the system $\{x_1, \ x_2, \ \ldots, \ x_j\}$ forms a basis of the $j$-dimensional subspace $L$. Examine the exterior product
$x_1 \wedge x_2 \wedge \ldots \wedge x_j$. Since $x_1, \ x_2, \ \ldots, \ x_j$ are linearly independent, the element $x_1 \wedge x_2 \wedge
\ldots \wedge x_j$ is nonzero and belongs to the line ${\mathcal
A}_j(L) \subset (K_j\cup (- K_j))$. Since $x_1 \wedge x_2 \wedge \ldots
\wedge x_j \in (K_j\cup (- K_j)) \setminus \{0\}$ for some nonzero
vectors $x_2, \ \ldots, \ x_j \in {\mathbb R}^n$,
we have $x_1 \in \widehat{T}(K_j)$.

$\Rightarrow$ The inclusion $0 \in L$ is obvious for any subspace $L \subset
{\mathbb R}^n$. Let $x_1 \in \widehat{T}(K_j)$ be nonzero. Then there exist
nonzero vectors $x_2, \ \ldots, \ x_j \in {\mathbb R}^n$ for
which $x_1 \wedge x_2 \wedge \ldots \wedge x_j \in (K_j \cup (-
K_j)) \setminus \{0\}$. Since $x_1 \wedge x_2 \wedge \ldots \wedge x_j \neq 0$, they are linearly independent. Examine the $j$-dimensional subspace $L =
{\rm Lin}(x_1, \ x_2, \ \ldots, \ x_j)$. Since $K_j \cup (-
K_j)$ is a cone of rank $1$ in $\wedge^j {\mathbb R}^n$, the line $\{t(x_1 \wedge x_2 \wedge \ldots \wedge x_j)\}_{t \in
{\mathbb R}}$ corresponding to the subspace $L$ belongs to
$K_j\cup (- K_j)$.

The second part of the lemma is proved analogically. \end{pf}

Now examine the set $\widetilde{T}(K_j)$ defined in the following way:

$$\widetilde{T}(K_j) = \{x_1 \in {\mathbb{R}}^{n}: \exists  \
x_2, \ \ldots, \ x_j \in {\mathbb{R}}^{n}, \ \mbox{for which}$$ $$ \
x_1 \wedge x_2 \wedge \ldots \wedge x_j \in ({\rm int}(K_j) \cup
{\rm int}(- K_j)) \}.$$

The above definition implies that $\widetilde{T}(K_j)
\subset T(K_j)$ and $\widetilde{T}(K_j)
\subset \widehat{T}(K_j)$.
The following statement describes the relations between the sets $T(K_j)$, $\widehat{T}(K_j)$
and $\widetilde{T}(K_j)$.

\begin{thm}
Let $K_j \subset \wedge^j {\mathbb{R}}^{n}$ be a proper cone. Then $$  \widetilde{T}(K_j) = {\rm int}(T(K_j)); \eqno(1)$$
$$T(K_j) \subseteq \widehat{T}(K_j). \eqno(2)$$
\end{thm}

\begin{pf}
 To prove (1), it is enough to show that the set $\widetilde{T}(K_j)$ is open.
Let $x_1 \in \widetilde{T}(K_j)$. Then there
exist elements $x_2, \ \ldots, \ x_j \in {\mathbb{R}}^{n}$ and
a number $r
> 0$ such that $B(x_1 \wedge x_2 \wedge \ldots \wedge x_j, r) \subset ({\rm
int}(K_j) \cup {\rm int}(- K_j))$. Let us find a number $r'
>0$ such that $B(x_1, r') \subset \widetilde{T}(K_j)$. Take $ r' =
\frac{r}{2j!\|x_2\| \ldots \|x_j\|}$. Indeed, the following
inequalities hold for every $x_1' \in B(x_1, r')$:
$$\|x_1' \wedge x_2 \wedge \ldots \wedge x_j - x_1 \wedge x_2 \wedge \ldots \wedge x_j\|
= \|(x_1' - x_1) \wedge x_2 \wedge \ldots \wedge x_j\| \leq $$
$$ \leq j! \|x_1' - x_1\|\|x_2\| \ldots \|x_j\| < j! \frac{r}{j!\|x_2\|
\ldots \|x_j\|}\|x_2\| \ldots \|x_j\| = r.$$

Since $x_1' \wedge x_2 \wedge \ldots \wedge x_j \in B(x_1
\wedge x_2 \wedge \ldots \wedge x_j, r) \subset ({\rm int}(K_j) \cup
{\rm int}(- K_j))$, we have $x_1' \in \widetilde{T}(K_j)$.

To prove (2), it is enough to show that the set $\widehat{T}(K_j)$ is closed. Let us take a sequence $\{x_1^n\}_{n=1}^{\infty} \in \widehat{T}(K_j)$ which converges to a nonzero element $x_1 \in {\mathbb R}^n$.
To show that $x_1 \in \widehat{T}(K_j)$, let us take
elements $x_2, \ \ldots, \ x_j \in {\mathbb R}^n$ for which $x_1
\wedge x_2 \wedge \ldots \wedge x_j \in (K_j \cup (- K_j))
\setminus \{0\}$.

Since the elements $x_1^n$ belong to $\widehat{T}(K_j)$ for every
$n = 1, \ 2, \ \ldots$, we can find elements $x_2^n, \
\ldots, \ x_j^n \in {\mathbb R}^n$ for which $x_1^n \wedge x_2^n
\wedge \ldots \wedge x_j^n \in (K_j \cup (- K_j)) \setminus
\{0\}$. Examine $j-1$ sequences $\{x_2^n\}_{n=1}^{\infty}, \
\ldots, \ \{x_j^n\}_{n=1}^{\infty}$. Without loss of
generality we can assume that $\|x_i^n\| = 1$ for
every $i = 2, \ \ldots, \ j$ and every $n = 1, \ 2, \ \ldots$.
Indeed, the linearity of the exterior product implies that
$$x_1^n \wedge \frac{x_2^n}{\|x_2^n\|}
\wedge \ldots \wedge \frac{x_j^n}{\|x_j^n\|} = \frac{1}{\|x_2^n\|
\ldots \|x_j^n\|}(x_1^n \wedge x_2^n \wedge \ldots \wedge x_j^n)
\in (K_j \cup (- K_j)) \setminus \{0\}.$$ Without loss of generality we can also assume that the linearly independent vectors $x_1^n, \ x_2^n, \
\ldots, \ x_j^n$ are mutually orthogonal for every $n = 1, \ 2, \ \ldots$.
Indeed, we can apply the Gram--Schmidt orthogonalization process to the set of $j$ linearly independent vectors $x_1^n, \ x_2^n, \
\ldots, \ x_j^n$. The obtained orthogonal vectors $\widetilde{x}_1^n, \ \widetilde{x}_2^n, \
\ldots, \ \widetilde{x}_j^n$ define the same $j$-dimensional subspace in
${\mathbb R}^n$ and the same line in $\wedge^j {\mathbb R}^n$.

Since all these sequences $\{x_2^n\}_{n=1}^{\infty}, \
\ldots, \ \{x_j^n\}_{n=1}^{\infty}$ are bounded, we can find a
converging subsequence in everyone of them. Let us take the necessary element $x_j$ equal to the limit of the corresponding converging subsequence
$x_j^{n_{k_j}}$. It is not difficult to see, that the elements $x_1, x_2, \ldots, x_j$ are nonzero and mutually orthogonal, so their exterior product $x_1 \wedge x_2 \wedge \ldots \wedge x_j$ is not equal to zero. Examine the sequence
$\{x_1^m \wedge x_2^m \wedge \ldots \wedge x_j^m\}_{m=1}^{\infty}$ of the exterior products of re-numbered elements of the subsequences $x_j^{n_{k_j}}$. The estimates $$\|x_1^m \wedge x_2^m \wedge \ldots \wedge x_j^m - x_1 \wedge x_2 \wedge \ldots \wedge x_j\| \leq $$ $$\leq  \sum_{k=1}^j\|x_1 \wedge \ldots \wedge x_{k-1} \wedge (x_k - x_k^m) \wedge x_{k+1}^m \ldots \wedge x_j^m \| \leq $$
$$ \leq \dfrac{1}{j!}\sum_{k=1}^j\|x_1\| \ldots  \|x_{k-1}\|\| x_k - x_k^m\| \|x_{k+1}^m \|\ldots \|x_j^m \| = $$ $$= \dfrac{1}{j!}(\| x_1 - x_1^m\| + \sum_{k=2}^j\|x_1\|\| x_k - x_k^m\|) $$ imply
that it converges to the element $x_1 \wedge x_2 \wedge \ldots \wedge x_j$. Since the sequence
$\{x_1^m \wedge x_2^m \wedge \ldots \wedge x_j^m\}_{m=1}^{\infty}$ belongs to the closed set $K_j\cup(-K_j)$ and the limit $x_1 \wedge x_2 \wedge \ldots \wedge x_j$ is nonzero, we conclude that $x_1 \wedge x_2 \wedge \ldots \wedge x_j \in (K_j \cup (- K_j)) \setminus \{0\}$. This implies that the element $x_1 \in \widehat{T}(K_j)$.
\end{pf}

\begin{thm} Let $K_j \subset \wedge^j{\mathbb R}^n$ be a proper cone. Let there exist a basis $e_1', \ \ldots, \ e_n'$ of ${\mathbb R}^n$ such that $K_j \subseteq K'_j$ where $K'_j \subset \wedge^j{\mathbb R}^n$ is one of the exterior basic cones defined by this
basis. Then the set $T(K_j)$, if it is not empty, is a cone of rank $j$.\end{thm}

\begin{pf} First let us prove that the set $T(K_j)$ is uniform, i.e. that for every $x_1 \in T(K_j)$,
$\alpha \in {\mathbb R}$ the element $\alpha x_1 \in T(K_j)$. It is enough to prove the above inclusion for every nonzero $\alpha$ and every $x_1 \in \widetilde{T}(K_j)$. Indeed, let
 $x_2, \ \ldots, \ x_j$ be nonzero elements for
which $x_1 \wedge x_2 \wedge \ldots \wedge x_j \in ({\rm int}(K_j) \cup
{\rm int}(- K_j))$. Let $\alpha$ be an arbitrary nonzero
number from ${\mathbb R}$. The linearity of the
exterior product implies that $\alpha x_1 \wedge x_2 \wedge \ldots \wedge
\frac{1}{\alpha}x_j = x_1 \wedge x_2 \wedge \ldots \wedge
x_j \in ({\rm int}(K_j) \cup
{\rm int}(- K_j))$.

 Lemma 3 implies that the set $T(K_j)$, if it is not empty, contains at least one $j$-dimensional
subspace. Then we have to prove, that there is no
$j+1$-dimensional subspace lying in $T(K_j)$. First
prove this fact for $n = j+1$. Let us show that $T(K_j)$
does not coincide with the whole of ${\mathbb R}^{j+1}$. Examine
the space ${\mathbb R}^{j+1}$ with the basis $e'_1, \ \ldots  , \
e'_{j+1}$ where $e'_1, \ \ldots , \ e'_{j+1}$ are given in the
condition of the theorem. Then the exterior products
$\{e'_{i_1}\wedge \ldots \wedge e'_{i_j}\}$ where $1 \leq i_1 < \ldots
< i_{j} \leq j+1$ form a basis in $\wedge^j {\mathbb R}^{j+1}$. So we
conclude that the cone $K'_j$ coincides with one of the basic cones of
the space $\wedge^j {\mathbb R}^{j+1}$. Without loss of generality we can assume that the maximal angle
$\theta_{\max}(K_j') = \frac{\pi}{2}$. Since $K_j \subseteq K_j'$,
we have $\theta_{\max}(K_j) \leq \theta_{\max}(K_j') =
\frac{\pi}{2}$. Let us examine the set ${\mathcal J}_{j+1}(K_j)$
where the operator ${\mathcal J}_{j+1}: \wedge^j {\mathbb R}^{j+1}
\rightarrow {\mathbb R}^{j+1}$ is defined in the following way:
$${\mathcal J}_{j+1} (e_{i_1}\wedge \ldots \wedge e_{i_{j}}) =
 (-1)^{k+1}e_k,$$
where $k = [j + 1]\setminus \{i_1, \ \ldots, \ i_{j}\}$.
Since the operator ${\mathcal J}_{j+1}$ is linear and invertible,
we conclude that ${\mathcal J}_{j+1}(K_j)$ is a proper cone in ${\mathbb
R}^{j+1}$. Moreover, ${\mathcal J}_{j+1}(K_j)$ belongs to
${\mathcal J}_{j+1}(K'_j)$ which coincides with one of the basic
cones of ${\mathbb R}^{j+1}$. So we conclude that
$\theta_{\max}({\mathcal J}_{j+1}(K_j)) \leq
\theta_{\max}({\mathcal J}_{j+1}(K'_j)) = \frac{\pi}{2}$. Lemma 3 and the properties of the operator ${\mathcal
J}_{j+1}$ imply that the set $T(K_j)$ is enclosed to the set
of all hyperplanes which orthogonal lines belong to
${\mathcal J}_{j+1}(K_j)\cup (- {\mathcal J}_{j+1}(K_j))$. Let us
show that ${\rm int}({\mathcal J}_{j+1}(K_j))\cup{\rm int}({-
\mathcal J}_{j+1}(K_j))$ does not belong to $T(K_j)$.
Indeed, let $x \in ({\rm int}({\mathcal J}_{j+1}(K_j))\cup{\rm
int}({- \mathcal J}_{j+1}(K_j))) \cap T(K_j)$. Since $x \in T(K_j)$, there exists a hyperplane $L$ such that $x \in L$ and the line $n$ orthogonal to $L$ belongs to ${\rm int}({\mathcal J}_{j+1}(K_j))\cup{\rm
int}({- \mathcal J}_{j+1}(K_j))$. Examine the angle $\theta$
between $n$ and $x$. It is equal to $\frac{\pi}{2}$. However, $x
\in {\rm int}({\mathcal J}_{j+1}(K_j))$, so the inequality $\theta
< \theta_{\max}({\mathcal J}_{j+1}(K_j)) \leq \frac{\pi}{2}$ holds. We came to the contradiction. Thus the set $T(K_j)$
does not coincide with the whole of ${\mathbb R}^{j+1}$. For $n =
j+1$ the theorem is proved.

Now let $n > j+1$. Let us prove the theorem by contradiction. Let $x_1, \
\ldots, \ x_{j+1}$ be $j+1$ linearly independent vectors, any
linear combination $c_1x_1 + \ldots + c_{j+1}x_{j+1} $ $ \ (c_1, \
\ldots, \ c_{j+1} \in {\mathbb R})$ of which belongs to $T(K_j)$. Let $(x_i^1, \ \ldots, x_i^n)$ be the coordinates of the
vector $x_i \ (i = 1, \ \ldots, \ j+1)$ in the basis $e'_1, \
\ldots, \ e'_n$. If the vectors $x_1, \ \ldots,
\ x_{j+1}$ are linearly independent, then at least one of the
minors of the form $
\begin{vmatrix}
  x_1^{i_1} & \ldots & x_{j+1}^{i_1} \\
  \ldots & \ldots & \ldots \\
  x_1^{i_{j+1}} & \ldots & x_{j+1}^{i_{j+1}} \\
\end{vmatrix}
$ where $1 \leq i_1 < \ldots < i_{j+1} \leq n$ is not equal to
zero. Examine a $j+1$-dimensional subspace $L$ of the space
${\mathbb R}^n$ spanned by the basic vectors $e'_{i_1}, \ \ldots,
\ e'_{i_{j+1}}$ and the corresponding subspace $\wedge^j L$ of
the space ${\mathbb R}^{\binom{n}j}$ spanned by all the possible
exterior products of the basic vectors $e'_{i_1}, \ \ldots, \
e'_{i_{j+1}}$. Examine a projection ${\rm pr}_{\wedge^j L} K_j$ of
the cone $K_j$ to the subspace $\wedge^j L$. It is not difficult
to see, that ${\rm pr}_{\wedge^j L} K_j \subseteq {\rm pr}_{
\wedge^j L} K'_j$ which is an exterior basic cone in $\wedge^j L$ (see Section 3, property 1 of exterior basic cones).
Since the space $L$ is $(j+1)$-dimensional, the statement of the
theorem holds. So we can find a vector $\varphi \in (L \setminus T({\rm pr}_{\wedge^j L} K_j))$ with
the coordinates $(\varphi_{i_1}, \ \ldots, \ \varphi_{i_{j+1}})$.
 Examine
the system
$$ \left\{\begin{array}{ccc}
 c_1x_1^{i_1} + \ldots + c_{j+1}x_{j+1}^{i_1} = \varphi_{i_1};
\\ \ldots \ \ \ldots \ \ \ldots
\\ c_1x_1^{i_{j+1}} + \ldots + c_{j+1}x_{j+1}^{i_{j+1}} = \varphi_{i_{j+1}} \end{array}\right.$$
This system has a unique solution $(c_1^0, \ \ldots, \
c_{j+1}^0)$. The vector $y_1 = c_1^0x_1 + \ldots + c_{j+1}^0x_{j+1}$
belongs to $T(K_j)$. Examine the case when we can find nonzero vectors
$y_2, \ \ldots, \ y_j \in {\mathbb R}^n$ such that $y_1 \wedge
y_2 \wedge \ldots \wedge y_j \in ({\rm int}(K_j) \cup
{\rm int}(- K_j))$. In this case ${\rm pr}_{\wedge^j L}(y_1 \wedge
y_2 \wedge \ldots \wedge y_j) \in {\rm
pr}_{\wedge^j L}({\rm int}(K_j) \cup
{\rm int}(- K_j))$ and is obviously nonzero. However, since ${\rm pr}_{\wedge^j L}(y_1 \wedge
y_2 \wedge \ldots \wedge y_j) = {\rm pr}_{L}(y_1) \wedge {\rm
pr}_{L}(y_2) \wedge \ldots \wedge {\rm pr}_{L}(y_j)$ and ${\rm pr}_{L}(y_1) = \varphi$, the above inclusion implies that $\varphi \in T({\rm pr}_{\wedge^j L}
K_j)$. In the case when the vector $y_1 = c_1^0x_1 + \ldots + c_{j+1}^0x_{j+1}$ is the limit of the converging sequence from $\widetilde{T}(K_j)$ we can construct the sequence of elements from $T({\rm pr}_{\wedge^j L}
K_j)$ converging to $\varphi$. So $\varphi \in T({\rm pr}_{\wedge^j L}
K_j)$. We came to the contradiction. \end{pf}

\begin{rmk}
 Note, that we do not use the convexity of the cone $K_j$ in the proof of Theorem 5.
\end{rmk}

\section{Construction of other cones of finite ranks}

Given $n$ proper cones $K_1 \subset {\mathbb R}^n, \ K_2 \subset
\wedge^2{\mathbb R}^n, \ \ldots, \ K_n \subset \wedge^n{\mathbb R}^n$.

Let us define successively the following sets $T(K_1, \ \ldots, \ K_j)$.

$$T(K_1) = K_1 \cup (- K_1);$$

$$T(K_1, K_2) = \overline{\{x_1 \in {\mathbb{R}}^{n}: \exists  \ x_2 \in
T(K_1),}$$ $$\overline{ \mbox{for which} \ x_1 \wedge x_2 \in
({\rm int}(K_2) \cup
{\rm int}(- K_2))\}};$$

$$T(K_1, K_2, K_3) = \overline{\{x_1 \in {\mathbb{R}}^{n}: \exists  \ x_2 \in
T(K_1), \ x_3 \in
T(K_1, K_2) }$$ $$ \overline{\mbox{for which} \ x_1 \wedge x_2 \wedge x_3 \in
({\rm int}(K_3) \cup
{\rm int}(- K_3))\}};$$

$$\ldots \ \ldots \ \ldots $$

$$T(K_1,\ldots, K_j) = \overline{\{x_1 \in {\mathbb{R}}^{n}: \exists  \ x_2 \in
T(K_1), \ \ldots, \ x_j \in
T(K_1, \ldots, K_{j-1})} $$ $$\overline{ \mbox{for which} \ x_1 \wedge x_2 \wedge \ldots \wedge x_j \in
({\rm int}(K_j) \cup
{\rm int}(- K_j))\}};$$

$$\ldots \ \ldots \ \ldots $$

$$T(K_1,\ldots, K_n) = \overline{\{x_1 \in {\mathbb{R}}^{n}: \exists  \ x_2 \in
T(K_1), \ \ldots, \ x_n \in
T(K_1, \ldots, K_{n-1})} $$ $$ \overline{\mbox{for which} \ x_1 \wedge x_2 \wedge \ldots \wedge x_n \in
({\rm int}(K_n) \cup
{\rm int}(- K_n))\}}.$$

The inclusions $T(K_1, \ldots, K_j) \subseteq T(K_j)$ are obvious for any proper cones $K_1 \subset {\mathbb R}^n, \ \ldots, \ K_j \subset \wedge^j{\mathbb R}^n$ and any $j = 2, \ \ldots, \ n$. The following theorem describes the structure of the sets $T(K_1, \ldots, K_j)$.

\begin{thm}
Let $K_1 \subset {\mathbb R}^n, \ K_2 \subset
\wedge^2{\mathbb R}^n, \ \ldots, \ K_n \subset \wedge^n{\mathbb R}^n$ be proper cones. Let for every $j \ \ (j = 2, \ \ldots, \
n)$ there exist a basis $e_1'(j), \ \ldots, \ e_n'(j)$ of ${\mathbb R}^n$ such that $K_j \subseteq K'_j$ where $K'_j \subset \wedge^j {\mathbb R}^n$ is one of the exterior basic cones defined by this basis. Then every set $T(K_1,\ldots, K_j)$, if it is not empty, is a cone of rank $j$.
\end{thm}

\begin{pf}
 First let us prove that for every $x_1 \in T(K_1, \ldots, K_j)$,
$\alpha \in {\mathbb R}$ the element $\alpha x_1 \in T(K_1, \ldots, K_j)$.
Let $x_1$ be an arbitrary element from $T(K_1, \ldots, K_j)$. Without loss of generality we can assume that there
exist nonzero elements  $x_2 \in
T(K_1), \ \ldots, \ x_j \in
T(K_1, \ldots, K_{j-1})$ for which $x_1
\wedge x_2\wedge \ldots \wedge x_j \in ({\rm int}(K_j) \cup
{\rm int}(- K_j))$. (Otherwise we shall consider $x_1$ as a limit of the converging sequence of the elements which satisfy the above condition.)  Let $\alpha$ be an arbitrary
nonzero number from ${\mathbb R}$. Since the set $T(K_1)$ is obviously uniform and $x_2 \in T(K_1)$, the element
$\frac{1}{\alpha}x_2$ also belongs to $T(K_1)$. The linearity of the exterior product implies that $\alpha x_1 \wedge
\frac{1}{\alpha}x_2 \wedge \ldots \wedge x_j = x_1
\wedge x_2\wedge \ldots \wedge x_j \in ({\rm int}(K_j) \cup
{\rm int}(- K_j))$. It is obvious that $0 \in T(K_1, \ldots, K_j)$ for $\alpha
= 0$. Hence $\alpha x_1 \in
T(K_1, \ldots, K_j)$ for every $\alpha \in {\mathbb R}$.

The definition of the set $T(K_1,\ldots, K_j)$ implies that $T(K_1,\ldots, K_j)$ is closed. The  inclusion $T(K_1, \ldots, K_j) \subseteq T(K_j)$ implies that $T(K_1,\ldots, K_j)$ does not contain any $j+1$-dimensional subspace. Let us show that the set $T(K_1, \ldots, K_j)$, if it is nonempty, contains at least one $j$-dimensional subspace. Indeed, let there exist at least one $x_1 \in {\mathbb R}^n$ such that we can find nonzero elements $x_2 \in
T(K_1), \ \ldots, \ x_j \in
T(K_1, \ldots, K_{j-1})$ for which $x_1
\wedge x_2\wedge \ldots \wedge x_j \in ({\rm int}(K_j) \cup
{\rm int}(- K_j))$. Let us prove that the $j$-dimensional subspace $L = {\rm Lin}(x_1, x_2, \ldots, x_j)$ belongs to $T(K_1, \ldots, K_j)$. Indeed,
examine the linear combination $\sum\limits_{i=1}^jc_ix_i$ where $c_1, \ \ldots, \ c_j \in {\mathbb{R}}, \ \sum\limits_{i=1}^jc_i^2 \neq 0$. If $c_1 \neq 0$, then we have the equality $(\sum\limits_{i=1}^jc_ix_i)\wedge x_2\wedge \ldots \wedge x_j = c_1(x_1
\wedge x_2 \wedge \ldots \wedge x_j) \in ({\rm int}(K_j) \cup
{\rm int}(- K_j))$. So $\sum\limits_{i=1}^jc_ix_i \in T(K_1, \ \ldots, \ K_j)$ for any $c_1 \neq 0$. If $c_1=0$, then there exists a sequence of the elements $\{\dfrac{1}{n}x_1 + \sum\limits_{i=2}^jc_ix_i\}_{n = 1}^\infty \in L$ which
converges to the vector $\sum\limits_{i=2}^jc_ix_i$. It follows from the above reasoning, that $\dfrac{1}{n}x_1 + \sum\limits_{i=2}^jc_ix_i \in
T(K_1, \ldots, K_j)$ for every $n$. So $\sum\limits_{i=2}^jc_ix_i \in T(K_1, \ldots, K_j)$.
\end{pf}

{\bf Example 1.} Let us examine the space ${\mathbb R}^3$ with a fixed basis $e_1, \ e_2, \ e_3$. Let $K_1 \subset {\mathbb R}^3$ be a basic cone spanned by the vectors $e_1, \ e_2, \ -e_3$. Let $K_2 \subset \wedge^2 {\mathbb R}^3$ be an exterior basic cone spanned by the exterior products $e_1 \wedge e_2, \ e_3
\wedge e_1, \ e_2 \wedge e_3$. In this
case the set $T(K_2)$ is a cone of rank $2$ which coincides
with the set ${\mathbb R}^3 \setminus ({\rm int}({\mathcal
J}(K_2))\cup {\rm int}(- {\mathcal J}(K_2)))$. The cone ${\mathcal
J}(K_2)$ is spanned by the vectors $e_1, \ e_2, \ e_3$, i.e. is
equal to ${\mathbb R}^3_+$. So the cone $T(K_1)$ coincides
with ${\mathbb R}^3 \setminus ({\rm int}({\mathbb R}^3_+)\cup {\rm
int}(- {\mathbb R}^3_+))$, i.e. with the set of all vectors which
have at least two coordinates of different signs or at least one
zero coordinate. It is not difficult to see, that the set $T(K_1,
K_2)$ is the set of all $2$-dimensional
subspaces $L \subset {\mathbb R}^3$ which satisfy the following
conditions.
\begin{enumerate}
\item[\rm 1.] The corresponding line ${\mathcal A}(L)$
belongs to $K_2\cup (- K_2)$;
\item[\rm 2.] The intersection $L \cap K_1 \neq \{0\}$.
\end{enumerate}
Since every $2$-dimensional subspace from $T(K_2)$ has a nonzero intersection with the cone $K_1$, we have the equality $T(K_1, K_2) = T(K_2)$.

{\bf Example 2.} Let $K_1$ be ${\mathbb R}^3_+$, $K_2 \subset \wedge^2 {\mathbb R}^3$ be an exterior basic cone spanned by the exterior products $e_1 \wedge e_2, \ e_3
\wedge e_1, \ e_2 \wedge e_3$. In this case, as it was shown above, $T(K_2) = {\mathbb R}^3 \setminus ({\rm int}({\mathbb R}^3_+)\cup {\rm
int}(- {\mathbb R}^3_+))$. It is easy to see, that the set $T(K_1,K_2) \subset T(K_2)$ is the union of three basic subspaces ${\rm Lin}(e_1,e_2)$, ${\rm Lin}(e_1,e_3)$ and ${\rm Lin}(e_2,e_3)$, i.e. the set of all vectors which have at least one zero coordinate. Note, that if $K_1$ is any proper cone which satisfies the inclusion $K_1 \subset {\rm int}({\mathbb R}^3_+)$, we obtain the equality $T(K_1,K_2) = \emptyset$.

Let us define successively the following sets $\widetilde{T}(K_1, \ \ldots, \ K_j)$.

$$\widetilde{T}(K_1) = {\rm int}(K_1) \cup (- {\rm int}(K_1));$$

$$\widetilde{T}(K_1, K_2) = \{x_1 \in {\mathbb{R}}^{n}: \exists  \ x_2 \in
\widetilde{T}(K_1),$$ $$ \ \mbox{for which} \ x_1 \wedge x_2 \in
({\rm int}(K_2) \cup
{\rm int}(- K_2))\};$$

$$\widetilde{T}(K_1, K_2, K_3) = \{x_1 \in {\mathbb{R}}^{n}: \exists  \ x_2 \in
\widetilde{T}(K_1), \ x_3 \in
\widetilde{T}(K_1, K_2) $$ $$ \mbox{for which} \ x_1 \wedge x_2 \wedge x_3 \in
({\rm int}(K_3) \cup
{\rm int}(- K_3))\};$$

$$\ldots \ \ldots \ \ldots $$

$$\widetilde{T}(K_1,\ldots, K_j) = \{x_1 \in {\mathbb{R}}^{n}: \exists  \ x_2 \in
\widetilde{T}(K_1), \ \ldots, \ x_j \in
\widetilde{T}(K_1, \ldots, K_{j-1}) $$ $$ \mbox{for which} \ x_1 \wedge x_2 \wedge \ldots \wedge x_j \in
({\rm int}(K_j) \cup
{\rm int}(- K_j))\};$$

$$\ldots \ \ldots \ \ldots $$

$$\widetilde{T}(K_1,\ldots, K_n) = \{x_1 \in {\mathbb{R}}^{n}: \exists  \ x_2 \in
\widetilde{T}(K_1), \ \ldots, \ x_n \in
\widetilde{T}(K_1, \ldots, K_{n-1}) $$ $$ \mbox{for which} \ x_1 \wedge x_2 \wedge \ldots \wedge x_n \in
({\rm int}(K_n) \cup
{\rm int}(- K_n))\}.$$

The inclusions $\widetilde{T}(K_1, \ldots, K_j) \subseteq \widetilde{T}(K_j)$ and $\widetilde{T}(K_1, \ldots, K_j) \subset T(K_1, \ldots, K_j)$ are obvious for any proper cones $K_1 \subset {\mathbb R}^n, \ \ldots, \ K_j \subset \wedge^j{\mathbb R}^n$ and any $j = 2, \ \ldots, \ n$. The following theorem describes the structure of the sets $\widetilde{T}(K_1, \ldots, K_j)$.

\begin{thm}
Let $K_1 \subset {\mathbb R}^n, \ K_2 \subset
\wedge^2{\mathbb R}^n, \ \ldots, \ K_n \subset \wedge^n{\mathbb R}^n$ be proper cones. Then the following inclusion is valid for every $j \ \ (j = 2, \ \ldots, \
n)$: $$\widetilde{T}(K_1,  \ldots,  K_j) \subseteq {\rm int}(T(K_1, \ldots, K_j)).$$
\end{thm}

\begin{pf}
It is enough for the proof to show that if the set $\widetilde{T}(K_1,  \ldots, K_j)$ is nonempty, then it is open. The proof is analogous to the proof of the first part of Theorem 4.
\end{pf}

Note, that it is also not difficult to define successively the sets $\widehat{T}(K_1, \ \ldots, \ K_j)$, $j = 2, \ \ldots, \ n$, using the formulae:

$$\widehat{T}(K_1) = K_1 \cup (- K_1);$$

$$\widehat{T}(K_1,\ldots, K_j) = $$ $$ = \overline{\{x_1 \in {\mathbb{R}}^{n}: \exists  \ x_2 \in
\widehat{T}(K_1)\setminus \{0\}, \ \ldots, \ x_j \in
\widehat{T}(K_1, \ldots, K_{j-1})\setminus \{0\}} $$ $$ \overline{\mbox{for which} \ x_1 \wedge x_2 \wedge \ldots \wedge x_j \in
(K_j\cup (- K_j)) \setminus \{0\}\}}.$$

\section{Cone-preserving maps in ${\mathbb{R}}^{n}$}
One of the most important results of the theory of nonnegative matrices is the famous Perron--Frobenius theorem. This theorem states the existence of the greatest in absolute value positive eigenvalue with the corresponding nonzero nonnegative eigenvector of a nonsingular nonnegative matrix (see, for example, \cite{BERPL}, p. 26).
Later we shall use the cone-theoretic generalizations of this result. So let us recall some definitions and statements of the theory of cone-preserving maps
(see \cite{BERPL, TAM}).

Let $K \subset {\mathbb R}^n$ be a proper cone. A linear operator
$A:{\mathbb{R}}^{n} \rightarrow {\mathbb{R}}^{n}$ is called {\it
$K$-positive} or {\it positive with respect to the cone $K$} if
$A(K \setminus \{0\}) \subseteq {\rm int}(K)$. In the case of $K =
{\mathbb{R}}^{n}_+$ $K$-positive operators are called {\it
positive}. It is easy to see, that the operator is positive if and only if its matrix is positive.

Let us state the following generalization of the Perron
theorem (see \cite{BERPL}, p. 13, Theorem 3.26). Recall that an {\it eigenfunctional} of the operator $A$ is defined as an eigenvector of the adjoint operator $A^*$.

\begin{thm}[Generalized Perron] Let a linear operator $A : {\mathbb R}^n
\rightarrow {\mathbb R}^n$ be positive with respect to a proper
cone $K \subset {\mathbb R}^n$. Then:
\begin{enumerate}
\item[\rm 1.] The spectral radius $\rho(A)$ is a simple positive
eigenvalue of the operator $A$ different in absolute value from the
remaining eigenvalues.
\item[\rm 2.] The eigenvector $x_1$ corresponding to
the eigenvalue $\lambda_1 = \rho(A)$ belongs to ${\rm int}(K)$.
\item[\rm 3.] The eigenfunctional $x_1^*$ corresponding to
the eigenvalue $\lambda_1 = \rho(A)$ belongs to ${\rm int}(K^*)$ (i.e. satisfies the inequality
$\langle x, x_1^* \rangle > 0$ for every nonzero $x \in K$).
\end{enumerate}
\end{thm}

A linear operator
$A:{\mathbb{R}}^{n} \rightarrow {\mathbb{R}}^{n}$ is called {\it
$K$-nonnegative} or {\it nonnegative with respect to the cone $K$} if it leaves the cone $K$ invariant (i.e.
$AK \subseteq K$). In the case of $K =
{\mathbb{R}}^{n}_+$ $K$-nonnegative operators are called {\it
nonnegative}.

Recall the following obvious fact.
\begin{lem} A linear operator
$A:{\mathbb{R}}^{n} \rightarrow {\mathbb{R}}^{n}$ is $K$-nonnegative ($K$-positive) if and only if the adjoint operator $A^*:({\mathbb{R}}^{n})' \rightarrow ({\mathbb{R}}^{n})'$ is $K^*$-nonnegative (respectively, $K^*$-positive).
\end{lem}

A weaker version of the generalized Perron theorem is correct for \linebreak
$K$-nonnegative operators (see \cite{BERPL}, p. 6, Theorem 3.2).

\begin{thm} Let a linear operator $A : {\mathbb R}^n
\rightarrow {\mathbb R}^n$ be nonnegative with respect to a proper
cone $K \subset {\mathbb R}^n$. Then:
\begin{enumerate}
\item[\rm 1.] The spectral radius $\rho(A)$ is a nonnegative
eigenvalue of the operator $A$.
\item[\rm 2.] The eigenvector $x_1$ corresponding to
the eigenvalue $\lambda_1 = \rho(A)$ belongs to $K$.
\item[\rm 3.] The eigenfunctional $x_1^*$ corresponding to
the eigenvalue $\lambda_1 = \rho(A)$ belongs to $K^*$.
\end{enumerate}
\end{thm}

Let us also state the ''inverse'' Perron theorem (see \cite{BERPL}, p. 8, Theorem 3.5 and p. 13, Theorem 3.26). Here $deg(\lambda)$ denotes the size of the largest diagonal block in the Jordan canonical
form of $A$ which contains $\lambda$.

\begin{thm}[Inverse Perron]
Let $\rho(A)$ be an eigenvalue of a linear operator $A : {\mathbb R}^n \rightarrow {\mathbb R}^n$. Let $deg (\lambda) \leq deg(\rho(A))$ for every eigenvalue $\lambda$ such that $|\lambda| = \rho(A)$. Then $A$ is nonnegative with respect to some proper cone $K$. Moreover, if $\rho(A)$ is a simple eigenvalue of $A$ greater in absolute value than the remaining eigenvalues, then $A$ is positive with respect to some proper cone $K$.
\end{thm}

Examine a subclass in the class of $K$-nonnegative operators which has the same spectral properties that $K$-positive operators. A linear operator $A: {\mathbb R}^n \rightarrow {\mathbb
R}^n$ is called {\it $K$--primitive} or {\it primitive with respect to the cone $K$}, if $AK \subseteq K$ and the only nonempty
subset of $\partial(K)$ which is left invariant by $A$ is $\{0\}$.
This definition was given by Barker (see \cite{BARK}, see also
\cite{TAM}). The following statement is correct (see \cite{BERPL}, p. 18, Corollary 4.13).
\begin{lem}
A linear operator $A: {\mathbb R}^n \rightarrow {\mathbb R}^n$ is primitive with respect to some proper cone $K$ if and only if there exists a proper cone $\widetilde{K}$ such that $A$ is positive with respect to $\widetilde{K}$.
\end{lem}

Let us examine the operators which leave invariant basic cones in ${\mathbb R}^n$. Every basic cone can be transformed into the cone
${\mathbb R}^n_+$ by a simple linear transformation with a diagonal transformation matrix. Thus {\it a linear operator $A: {\mathbb R}^n \rightarrow {\mathbb R}^n$ leaves invariant a basic cone in ${\mathbb R}^n$ if and only if its matrix ${\mathbf A}$ can be represented in the
following form:
$${\mathbf A} = {\mathbf D} \widetilde{{\mathbf A}} {\mathbf
D}^{-1}, $$ where $\widetilde{{\mathbf A}}$ is a
nonnegative matrix, ${\mathbf D}$ is a diagonal matrix, which
diagonal elements are equal to $\pm 1$.}

 Let $J$ be any subset of $[n] := \{1, 2, \ldots, n\}$. Then $J^c:=[n] \setminus J$ and
$$[n] \times [n] = (J\times J)\cup(J^c \times J^c)\cup (J \times J^c) \cup (J^c \times J)$$ is a partition of $[n] \times [n]$ into four pairwise disjoint subsets.

A matrix ${\mathbf A} = \{a_{ij}\}_{i,j=1}^n$ is called {\it J-sign-symmetric (JS)} if
$$a_{ij} \geq 0 \quad \mbox{on} \quad (J \times J)\cup (J^c \times J^c);$$
and
$$a_{ij} \leq 0 \quad  \mbox{on} \quad (J \times J^c)\cup (J^c \times J).$$

 A matrix ${\mathbf A} = \{a_{ij}\}_{i,j=1}^n$ is called {\it strictly J-sign-symmetric (SJS)} if
$$a_{ij} > 0 \quad \mbox{on} \quad (J \times J)\cup (J^c \times J^c);$$
and
$$a_{ij} < 0 \quad  \mbox{on} \quad (J \times J^c)\cup (J^c \times J).$$

It is easy to see, that the number of all different types of strictly J--sign-symmetric $n \times n$ matrices is
equal to the number of basic cones divided by $2$, i.e. $2^{n - 1}$.

We recall a simple fact that a matrix is diagonally similar to a nonnegative (positive) matrix if and only if it is J-sign-symmetric (respectively, strictly J-sign-symmetric) (see, e.g., \cite{KU1}). Thus a linear operator $A: {\mathbb R}^n \rightarrow {\mathbb R}^n$ is nonnegative (positive) with respect to some basic cone $K$ if and only if its matrix $\mathbf A$ is J--sign-symmetric (respectively, strictly J--sign-symmetric).

\section{Exterior powers of operators in ${\mathbb R}^n$}

 Let us recall the following definitions and statements.

Let $A$ be a linear operator acting in the space ${\mathbb R}^n$. Then a linear operator $\wedge^j A$ acting in the space $\wedge^j
{\mathbb R}^n$ according to the rule:
$$ (\wedge^j A)(x_1 \wedge \ldots \wedge x_j ) = Ax_1 \wedge \ldots \wedge Ax_j,$$ is called
the {\it $j$th exterior power} of the operator $A$.

Later we shall use the following properties of $\wedge^j A$ (see,
for example, \cite{YOK}).
\begin{enumerate}
 \item[\rm \bf 1.] $\wedge^{j}A = 0$ if and only if $j \geq r$ where $r$ is the rank of the operator $A$.
\item[\rm \bf 2.] $\wedge^{j}I_{{\mathbb R}^n} = I_{\wedge^{j}{\mathbb R}^n}$ where $I$ is the identity operator.
\item[\rm \bf 3.] Let $A, B : {\mathbb R}^n \rightarrow {\mathbb R}^n$ be two linear operators. Then $\wedge
^{j}(A B) = (\wedge ^{j}A) (\wedge ^{j}B)$ (the Cauchy--Binet formula).
\item[\rm \bf 4.] The following equality is correct for every natural number $m$: $(\wedge ^{j}A)^{m} =
\wedge^{j}(A^{m}).$
\item[\rm \bf 5.] The $j$-th exterior power of an invertible operator is invertible and the following equality is correct: $(\wedge ^{j}A)^{- 1} =
\wedge^{j}(A^{- 1})$.
\item[\rm \bf 6.] Since $(\wedge^j{\mathbb R}^n)'$ can be considered as $\wedge^j({\mathbb R}^n)'$, the following equality for adjoint operators is correct: $(\wedge ^{j}A)^* = \wedge^{j}(A^*)$ (see \cite{YOK}, p. 89).
\end{enumerate}

If the operator $A$ is defined by the matrix ${\mathbf A} =
\{a_{ij}\}_{i,j = 1}^n$ in the basis $e_1, \ \ldots, \ e_n$, then
the matrix of $\wedge^j A$ in the basis
$\{e_{i_1} \wedge \ldots \wedge e_{i_j}\}$ where $1 \leq i_1<
\ldots < i_j \leq n$ coincides with the $j$th compound matrix $
{\mathbf A}^{(j)}$ of the initial matrix ${\mathbf A}$ (see, for example,
\cite{PINK}).

Recall the following statement concerning the eigenvalues of $\wedge^j A$.

\begin{thm}[Kronecker] Let $\{\lambda_{i}\}_{i = 1}^n$ be the set of
all eigenvalues of the operator $A: {\mathbb R}^n \rightarrow
{\mathbb R}^n$ repeated according to multiplicity. Then all the
possible products of the form $\{\lambda_{i_1} \ldots \lambda_{i_j} \}$, where $1
\leq i_1 < \ldots < i_j \leq n$ forms the set of all the possible eigenvalues
of the $j$th exterior power $\wedge^j A$ of the operator $A$ repeated
according to multiplicity.
\end{thm}

The Kronecker theorem is stated in terms of compound matrices
and proved without using exterior products in \cite{GANT} (see
\cite{GANT}, p. 80, Theorem 23).

\section{Generalized totally positive operators}

Let us give the basic definition of a generalized totally positive operator.
Let us fix $n$ proper cones $K_1 \subset {\mathbb R}^n, \ K_2 \subset \wedge^2 {\mathbb R}^n, \  \ldots, \ K_n \subset \wedge^n {\mathbb R}^n$.
Note, that the idea of fixing cones in the exterior powers of the initial space was first given by Yudovich in \cite{YD}. Slightly changing the terminology of \cite{YD}, we call the family of proper cones $K_1 \subset {\mathbb R}^n, \ K_2 \subset \wedge^2 {\mathbb R}^n, \  \ldots, \ K_n \subset \wedge^n {\mathbb R}^n$ a {\it totally positive structure} on the space ${\mathbb R}^n$.

A linear operator $A$ is called {\it generalized totally positive (GTP)} with respect to a totally positive structure $\{K_1, \ \ldots, \ K_n\}$ if
it is nonnegative with respect to the proper cone $K_1 \subset
{\mathbb R}^n$ and its $j$-th
exterior power $ \wedge^j A$ is nonnegative with respect to the
proper cone $K_j \subset \wedge^j {\mathbb R}^n$ for every $j$ $(j = 2, \ \ldots, \ n)$.

A linear operator $A$ is called {\it generalized strictly totally positive (GSTP)} with respect to a totally positive structure $\{K_1, \ \ldots, \ K_n\}$ if
it is positive with respect to the proper cone $K_1 \subset
{\mathbb R}^n$ and its $j$-th
exterior power $ \wedge^j A$ is positive with respect to the
proper cone $K_j \subset \wedge^j {\mathbb R}^n$ for every $j$ $(j = 2, \ \ldots, \ n)$.

In the case when $K_j = \wedge^j{\mathbb R}^n_+$ for every $j = 1, \ \ldots, \ n$, the definitions given above coincide with the
classical definitions of totally positive and strictly totally positive operators
given by Gantmacher and Krein in \cite{GANT}.

It is easy to see, that there exists at least one nonsingular GTP operator for every totally positive structure on the space ${\mathbb R}^n$ (for example, the identity operator $I$ which exterior powers according to Property 2 are also the identity operators). We call a totally positive structure {\it strictly totally positive} if there exists at least one strictly totally positive with respect to this structure operator.
Later we are going to show, that not every totally positive structure on ${\mathbb R}^n$ is strictly totally positive.

Now it is also possible to give the definition of generalized oscillatory operator which extends the class of oscillatory operators introduced by Gantmacher and Krein in \cite{GANT}.

A linear operator $A$ is called {\it generalized oscillatory (GO)} with respect to a totally positive structure $\{K_1, \ \ldots, \ K_n\}$ if
it is primitive with respect to the proper cone $K_1 \subset
{\mathbb R}^n$ and its $j$-th
exterior power $ \wedge^j A$ is primitive with respect to the
proper cone $K_j \subset \wedge^j{\mathbb R}^n$ for every $j$ $(j = 2, \ \ldots, \ n)$.

A linear operator $A$ is called {\it generalized sign-regular (GSR)} with respect to a totally positive structure $\{K_1, \ \ldots, \ K_n\}$ if there exist numbers
$\epsilon_1, \ \ldots, \ \epsilon_n$ each equal to $\pm 1$ such that
$\epsilon_1 A$ is nonnegative with respect to the proper cone $K_1 \subset
{\mathbb R}^n$ and $\epsilon_j \wedge^j A$ is nonnegative with respect to the
proper cone $K_j \subset \wedge^j {\mathbb R}^n$ for every $j$ $(j = 2, \ \ldots, \ n)$.

A linear operator $A$ is called {\it generalized strictly sign-regular (GSSR)} with respect to a totally positive structure $\{K_1, \ \ldots, \ K_n\}$ if there exists numbers
$\epsilon_1, \ \ldots, \ \epsilon_n$ each equal to $\pm 1$ such that
$\epsilon_1 A$ is positive with respect to the proper cone $K_1 \subset
{\mathbb R}^n$ and $\epsilon_j \wedge^j A$ is positive with respect to the
proper cone $K_j \subset \wedge^j {\mathbb R}^n$ for every $j$ $(j = 2, \ \ldots, \ n)$.

Let us fix a natural number $k$, $1 \leq k \leq n$ and choose
$k$ proper cones $K_1 \subset {\mathbb R}^n, \ K_2 \subset \wedge^2 {\mathbb R}^n, \  \ldots, \ K_k \subset \wedge^k {\mathbb R}^n$. We call the sequence of proper cones $K_1 \subset {\mathbb R}^n, \ K_2 \subset \wedge^2 {\mathbb R}^n, \  \ldots, \ K_k \subset \wedge^k {\mathbb R}^n$ a {\it $k$-totally positive structure} on the space ${\mathbb R}^n$.

A linear operator $A$ is called {\it generalized $k$-totally positive} with respect to a $k$-totally positive structure $\{K_1, \ \ldots, \ K_k\}$ if
it is nonnegative with respect to the proper cone $K_1 \subset
{\mathbb R}^n$ and its $j$-th
exterior power $ \wedge^j A$ is nonnegative with respect to the
proper cone $K_j \subset \wedge^j {\mathbb R}^n$ for every $j$ $(j = 2, \ \ldots, \ k)$.

We can easily give analogical definitions of generalized strictly $k$-totally positive, $k$-sign-regular and strictly $k$-sign-regular operators.

\section{Basic properties of GTP and GSTP operators}
Let us list some basic properties of GTP and GSTP operators.
\begin{prop}
Let a linear operator $A: {\mathbb R}^n \rightarrow {\mathbb R}^n$ be GTP (GSTP) with respect to a totally positive structure $\{K_1, \ \ldots, \ K_n\}$. Then $A^*$ (the ajoint of $A$) is GTP (respectively, GSTP) with respect to the totally positive structure $\{K_1^*, \ \ldots, \ K_n^*\}$.
\end{prop}
\begin{pf}
The proof follows from Lemma 9 and Property 6 of exterior powers (see Section 7).
\end{pf}
\begin{prop}
Let linear operators $A, B: {\mathbb R}^n \rightarrow {\mathbb R}^n$ be GTP with respect to a totally positive structure $\{K_1, \ \ldots, \ K_n\}$. Then $AB$ is also GTP with respect to the structure $\{K_1, \ \ldots, \ K_n\}$. If in this case one of the operators $A$ or $B$ is GSTP, while the other is nonsingular GTP, then $AB$ is GSTP with respect to the structure $\{K_1, \ \ldots, \ K_n\}$. In particular, if $A$ is GTP (GSTP) with respect to a totally positive structure $\{K_1, \ \ldots, \ K_n\}$, then the operator $A^m$ is GTP (respectively, GSTP) with respect to the same structure for every natural number $m$.
\end{prop}
\begin{pf}
 The proof follows from the Cauchy--Binet formula (see Section 7, Property 3 of exterior powers).
\end{pf}

 If the operators $A$ and $B$ are GTP with respect to different totally positive structures, then the statement of Proposition 15 may not be correct.

\begin{prop}
Let a linear operator $A: {\mathbb R}^n \rightarrow {\mathbb R}^n$ be GO with respect to a totally positive structure $\{K_1, \ \ldots, \ K_n\}$. Then $A$ is GSTP with respect to some other totally positive structure $\{\widetilde{K_1}, \ \ldots, \ \widetilde{K_n}\}$.
\end{prop}

\begin{pf}
 The proof follows from Lemma 12.
\end{pf}

Proposition 16 reduces the study of GO operators to the study of GSTP operators.

Now let us prove some propositions which describe the structure of the class of GTP operators.

\begin{prop}
Every linear operator $A: {\mathbb R}^n \rightarrow {\mathbb R}^n$ similar to a GTP (GSTP) operator is GTP (respectively, GSTP). In particular, if $A$ is similar to a TP (STP) operator, then $A$ is GTP (respectively, GSTP).
\end{prop}

\begin{pf} Let us represent the operator $A$ in the form $A = T\widetilde{A}T^{-1}$, where $T$ is nonsingular, $\widetilde{A}$ is GTP (GSTP) with respect to a totally positive structure $\{K_1, \ \ldots, \ K_n\}$. Then the Cauchy--Binet formula and Property 5 of the exterior powers imply that $A$ is GTP (GSTP) with respect to the totally positive structure $\{T(K_1), \ \ldots, \ (\wedge^nT) (K_n)\}$.
\end{pf}

Note, that the "inverse`` statement that every GTP (GSTP) operator is similar to some TP (STP) operator is not correct.

Proposition 17 shows that if all the eigenvalues of $A$ are positive and simple, then $A$ is GTP.
Indeed, since the Jordan canonical form $J$ of the operator $A$ is represented by a nonnegative diagonal matrix, we conclude that $J$ is TP. Then the equality $A = UJU^{-1}$ implies that $A$ is GTP with respect to the totally positive structure $\{U({\mathbb R}^n_+),  \ \ldots, \ (\wedge^nU) (\wedge^n{\mathbb R}^n_+)\}$.

Using the ''inverse`` Perron theorem (see Theorem 11) stated above, we can prove more general statements about operators with real spectrum. In this case the Jordan canonical form may not be totally positive or sign-regular.

\begin{prop}
Let all the eigenvalues of a linear operator $A: {\mathbb R}^n \rightarrow {\mathbb R}^n$ be real. Then $A$ is GSR. Moreover, if all the eigenvalues of $A$ are real, simple and different in absolute value from each other then $A$ is GSSR.
\end{prop}

\begin{pf}
Let us enumerate the eigenvalues of the operator $A$
in descending order of their absolute values (without taking into account their
multiplicities):
$$|\lambda_{1}| \geq | \lambda_{2}| \geq |\lambda_{3}| \geq \ldots
\geq |\lambda_{m}|.$$

Examine the greatest in absolute value eigenvalue $\lambda_{1}$. The reality of the spectrum implies that if the equality $|\lambda_{1}| = | \lambda_{2}|$ is true, then $\lambda_1 = - \lambda_2$. Assume that $deg(\lambda_1) \leq deg(\lambda_2)$, otherwise we shall re-number them. If the eigenvalue $\lambda_1$ is nonnegative, then it satisfies the conditions of Theorem 11. Applying Theorem 11 to the operator $A$, we obtain that $A$ is nonnegative with respect to some proper cone $K_1 \subset {\mathbb R}^n$. If the eigenvalue $\lambda_1$ is non-positive, then $- \lambda_1$ is nonnegative. Considering $- \lambda_1$ as the greatest in absolute value eigenvalue of the operator $-A$ we obtain that it satisfies the conditions of Theorem 11. Applying Theorem 11 to the operator $- A$ we obtain that $- A$ is nonnegative with respect to some proper cone $K_1 \subset {\mathbb R}^n$.

Examine the second exterior power $\wedge^2 A$. The Kronecker theorem implies $\wedge^2 A$ has no other eigenvalues, except all the possible products of the form $\lambda_{i_1}\lambda_{i_2}$, where $1 \leq i_1 < i_2 \leq n$. Therefore the greatest in absolute value eigenvalue $\lambda_1^{(2)}$ of $\wedge^2 A$ can be represented in the form of the product $\lambda_{i_1}\lambda_{i_2}$ with some values of
the indices $i_1,i_2$, \ $i_1 \leq i_2$. This representation implies that $\lambda_1^{(2)}$ is also real. Without loss of generality we can assume that $deg(\lambda_1^{(2)}) \leq deg(\lambda_2^{(2)})$ where $\lambda_2^{(2)}$ is any other eigenvalue of $\wedge^2 A$ equal in absolute value to $\lambda_1^{(2)}$. If $\lambda_1^{(2)}$ is nonnegative, we apply Theorem 11 to the operator $\wedge^2 A$, otherwise we apply Theorem 11 to the operator $- \wedge^2 A$. Thus we obtain that either $\wedge^2 A$ or $- \wedge^2 A$ is nonnegative with respect to some proper cone $K_2 \subset \wedge^2 {\mathbb R}^n$.

Repeating the above reasoning for $\wedge^j A$, $j = 3, \ \ldots, \ n$ we obtain that either the operator $\wedge^j A$ or $- \wedge^j A$ is nonnegative with respect to some proper cone $K_j \subset \wedge^j {\mathbb R}^n$.
So we have constructed the totally positive structure $\{K_1, \ \ldots, \ K_n\}$ such that the operator $A$ is GSR with respect to this structure.
The second part of the proposition is proved analogically.
\end{pf}

\begin{prop}
Let all the eigenvalues of a linear operator $A: {\mathbb R}^n \rightarrow {\mathbb R}^n$ be nonnegative. Then $A$ is TP. Moreover, if all the eigenvalues of $A$ are positive and simple then $A$ is GSTP.
\end{prop}
\begin{pf}
The proof follows from Proposition 18.
\end{pf}

Later we shall show that a GSTP operator always has a simple positive spectrum. However, the spectrum of a GTP operator may not be real.

Proposition 18 shows that the introduced above class of GSR operators covers the entire class of operators with real spectrum. This also implies that any operation which preserves the reality of the spectrum of an operator, preserves the class of GSR operators.

\begin{prop}
 Let a linear operator $A:{\mathbb R}^n \rightarrow {\mathbb R}^n$ be nonsingular GSR (GSSR). Then $A^{-1}$ is GSR (respectively, GSSR). In particular, if $A$ is nonsingular GTP (GSTP). Then $A^{-1}$ is GTP (respectively, GSTP).
\end{prop}
\begin{pf}
 The proof follows from Proposition 18 and Proposition 19.
\end{pf}

\section{Variation diminishing property of GTP operators}

Now we shall prove the generalization of the results by
Schoenberg concerning variation diminishing property of totally
positive matrices (see \cite{SCH1, SCH2}).
\begin{thm}
Let a linear operator $A: {\mathbb{R}}^{n} \rightarrow
{\mathbb{R}}^{n}$ be nonsingular GSR with respect to a totally positive structure $\{K_1, \ \ldots, \ K_n\}$. Then the following inclusions hold for every $j = 1, \ldots, \ n$: $$A(\widehat{T}(K_j)) \subseteq \widehat{T}(K_j); \eqno(3)$$
 $$A(\widehat{T}(K_1, \ldots, K_j)) \subseteq \widehat{T}(K_1, \ldots, K_j). \eqno(4)$$
\end{thm}

\begin{pf}
Let us assume that all the sets $\widehat{T}(K_j)$ are nonzero and all the sets $\widehat{T}(K_1, \ldots, K_j)$ are nonempty, otherwise the corresponding inclusions will be obvious.

Let $x_1$ be an arbitrary vector from $\widehat{T}(K_j)$. We prove that
$Ax_1 \in \widehat{T}(K_j)$. Let us find the elements $x_2, \ \ldots, \ x_j \in {\mathbb{R}}^{n}$ for which $x_1 \wedge x_2 \wedge \ldots \wedge x_j \in (K_j\cup (-K_j)) \setminus \{0\}$. Examine the elements $Ax_2, \ \ldots, \ Ax_j$. Since the
operator $\wedge^j A$ is nonsingular
nonnegative with respect to the cone $K_j$, we have $Ax_1 \wedge Ax_2 \wedge \ldots \wedge Ax_j = (\wedge^j A)(x_1 \wedge x_2 \wedge \ldots \wedge x_j) \in (K_j \cup (-K_j)) \setminus \{0\}$. Thus $Ax_1 \in \widehat{T}(K_j)$.

Let us prove inclusion (4) using the induction on $j$. The inclusion \linebreak $A(\widehat{T}(K_1)) \subseteq \widehat{T}(K_1)$ is obvious. Let us take $j = 2$ and prove that \linebreak
$A(\widehat{T}(K_1, K_2)) \subseteq \widehat{T}(K_1,K_2)$.
 Let there exists a nonzero element $x_2 \in \widehat{T}(K_1)$ such that $x_1 \wedge x_2 \in (K_2\cup (- K_2)) \setminus \{0\}$. Examine the element $Ax_2$. Since the
operator $A$ is nonsingular nonnegative with respect to the cone $K_1$,
we have $Ax_2 \in \widehat{T}(K_1)\setminus \{0\}$. Examine the
element $Ax_1 \wedge Ax_2$. Since the element $x_1 \wedge x_2 \in (K_2\cup
(- K_2)) \setminus \{0\}$ and the operator $\wedge^2 A$ is
nonsingular nonnegative with respect to the cone $K_2$, we have $Ax_1 \wedge
Ax_2 = (\wedge^2 A)(x_1 \wedge x_2) \in (K_2\cup (-K_2)) \setminus \{0\}$. Thus $Ax_1 \in \widehat{T}(K_1, K_2)$. Now let us consider the case when $x_1$ is the limit of a sequence
$\{x_1^n\}_{n = 1}^\infty$ such that there exists a nonzero element $x_2^n \in \widehat{T}(K_1)$ satisfying $ x^n_1
\wedge x^n_2 \in (K_2\cup (- K_2)) \setminus \{0\}\}$ for every element $x_1^n$. It follows from the above reasoning, that $Ax_1^n \in \widehat{T}(K_1,
K_2)$ for every $n = 1, \ 2, \ \ldots$. Since the sequence
$\{Ax_1^n\}_{n = 1}^ \infty$ converges to the vector $Ax_1$, we
conclude $Ax_1 \in \widehat{T}(K_1, K_2)$.

Let the statement of the theorem holds for $j-1$. Now let us prove inclusions (3) and (4) for $j$.
Let $x_1$ be an arbitrary vector from $\widehat{T}(K_1, \ldots, K_j)$. Prove, that
$Ax_1 \in \widehat{T}(K_1, \ldots, K_j)$. As it was shown above, without loss of generality we can assume that there
exist nonzero elements  $x_2 \in
\widehat{T}(K_1), \ \ldots,$ $x_j \in
\widehat{T}(K_1, \ldots, K_{j-1})$ such that $x_1
\wedge x_2\wedge \ldots \wedge x_j \in (K_j \cup (- K_j))\setminus \{0\}$.
Examine the elements $Ax_2, \ \ldots, \ Ax_j$ which are also nonzero. Using the inductive hypothesis, we obtain that $Ax_2 \in
\widehat{T}(K_1)\setminus\{0\}, \ \ldots, \ Ax_j \in
\widehat{T}(K_1, \ldots, K_{j-1})\setminus\{0\}$. Since the element $x_1 \wedge x_2 \wedge \ldots \wedge x_j \in (K_j\cup
(- K_j)) \setminus \{0\}$ and the operator $\wedge^j A$ is
nonsingular nonnegative with respect to the cone $K_j$, we have $Ax_1 \wedge
Ax_2 \wedge \ldots \wedge Ax_j = (\wedge^j A)(x_1 \wedge x_2 \wedge \ldots \wedge x_j) \in (K_j\cup (-K_j)) \setminus \{0\}$. Thus $Ax_1 \in \widehat{T}(K_1, \ldots K_j)$.
\end{pf}

\section{Gantmacher--Krein theorem for GSTP operators}

Let us state and prove the main theorem concerning spectral
properties of GSTP operators.

\begin{thm} Let a linear operator $A: {\mathbb{R}}^{n}
\rightarrow {\mathbb{R}}^{n}$ be GSTP with respect to a totally positive structure $\{K_1, \ \ldots, \ K_n\}$. Then all the eigenvalues of the operator $A$ are
positive and simple:
$$\rho(A)= \lambda_1 > \lambda_2 > \ldots > \lambda_n > 0.$$
 The first eigenvector $x_1$ corresponding to the maximal
eigenvalue $\lambda_1$ belongs to ${\rm int} (K_1)$ and the
$j$th eigenvector $x_j$ corresponding to the $j$th in absolute value
eigenvalue $\lambda_j$ belongs to ${\rm int}(T(K_1, \ldots,
K_j))\setminus T(K_1, \ldots, K_{j-1})$. Moreover, the following
inclusions hold:
$$\sum_{i=q}^p c_i x_i \in {\rm int}(T(K_1, \ \ldots, \
K_p))\setminus T(K_1, \ \ldots, \ K_{q - 1})$$ for each $1 \leq q \leq p
\leq n$ and $c_p \neq 0$;
$$\sum_{i=q}^p c_i x_i \in \overline{T(K_1, \ \ldots, \
K_p)\setminus T(K_1, \ \ldots, \ K_{q - 1})}$$ for each $1 \leq q \leq p
\leq n$.
\end{thm}

\begin{pf} The first part of the proof literally repeats the arguments used originally by Gantmacher and Krein. Let us list the eigenvalues of the operator $A$
in descending order of their absolute values (taking into account their
multiplicities):
$$|\lambda_{1}| \geq | \lambda_{2}| \geq |\lambda_{3}| \geq \ldots
\geq |\lambda_{n}|.$$

 Applying the generalized Perron theorem to the operator $A$ (which is positive with respect to the proper cone $K_1$),
we get: $\lambda_{1} = \rho(A)>0$ is a simple positive eigenvalue
of $A$, different in absolute value from the remaining eigenvalues.
The corresponding eigenvector $x_1$ belongs to ${\rm int}(K_1)$.
 Examine the second exterior power $\wedge^2 A$ which is positive with respect to the proper cone $K_2 \subset \wedge^2 {\mathbb R}^n$.
Applying generalized Perron theorem to $\wedge^2 A$, we get: $\rho(\wedge^2 A) > 0$ is
 a simple positive eigenvalue of $\wedge^2 A$, different in absolute value from the remaining eigenvalues. The corresponding eigenvector $\varphi_2$ belongs to ${\rm int}(K_2)$.

 As it follows from the statement of the Kronecker theorem, $\wedge^2 A$ has no other eigenvalues, except all the possible products of the form $\lambda_{i_1}\lambda_{i_2}$ where $1 \leq i_1 < i_2 \leq n$. Therefore $\rho(\wedge^2 A)>0$ can be represented in the
form of the product $\lambda_{i_1}\lambda_{i_2}$ with some values of
the indices $i_1,i_2$, \ $i_1 < i_2$. The facts that the
eigenvalues are listed in a descending order and there is only
one eigenvalue on the spectral circle $|\lambda| = \rho(A)$ imply that
 $\rho(\wedge^2 A) =  \lambda_{1}\lambda_{2} = \rho(A)\lambda_{2}$.
  Therefore $\lambda_{2} = \frac{\rho(\wedge^2 A)}{\rho(A)}>0$.

Repeating the same reasoning for $\wedge^j A$, $j = 3, \ \ldots, \ n$, we obtain the relations:
$$\lambda_{j} = \frac{\rho(\wedge^j A)}{\rho(\wedge^{j-1} A)}>0,$$
where $j = 3, \ \ldots, \ n$. The simplicity of the eigenvalues $\lambda_j$ for every $j$ also follows from the above relations and the simplicity of $\rho(\wedge^j A)$. Note, that the eigenvector $\varphi_j$ of the operator $\wedge^j A$ corresponding to the eigenvalue $\rho(\wedge^j A)$ belongs to ${\rm int}(K_j)$.

Now let us prove that the
$j$th eigenvector $x_j$ corresponding to the $j$th in absolute value
eigenvalue $\lambda_j$ belongs to ${\rm int}(T(K_1, \ \ldots, \
K_j))\setminus T(K_1, \ \ldots, \ K_{j-1})$. Let us prove this statement by induction on $j$. First take $j = 2$.

Since $\rho(\wedge^2 A)= \lambda_{1}\lambda_{2}$, the eigenvector
$\varphi_2$ corresponding to $\rho(\wedge^2 A)$ can be represented
in the form of the exterior product $\varphi_2 = x_1 \wedge x_2$ of
the first eigenvector $x_1 \in {\rm int}(K_1)$ and the second
eigenvector $x_2$. The inclusion $\varphi_2 \in {\rm int}(K_2)$ implies
$x_2 \in \widetilde {T}(K_1, K_2) \subseteq {\rm int}(T(K_1,
K_2))$.

 Let us show, that the vector
$x_2$ does not belong to $K_1 \cup (- K_1)$. It is enough for this
to show that $x_2$ belongs to the subspace $$X_2' = \{ x \in
{\mathbb R}^n: x_1^*(x) = 0 \},$$ where $x_1^*$ is the first
eigenfunctional of the operator $A$ corresponding to the maximal
eigenvalue $\lambda_1 = \rho(A)$. Indeed, $\langle x_1^*, x_2
\rangle = \frac{1}{\lambda_2}(\langle x_1^*, \lambda_2 x_2
\rangle) = \frac{1}{\lambda_2}(\langle x_1^*, Ax_2 \rangle) =
\frac{1}{\lambda_2}(\langle A^*x_1^*, x_2 \rangle) =
\frac{1}{\lambda_2}(\langle \lambda_1 x_1^*,  x_2 \rangle) =
\frac{\lambda_1}{\lambda_2}(\langle  x_1^*,  x_2 \rangle).$ Since
$\frac{\lambda_1}{\lambda_2} \neq 1$, the equality above is
valid if and only if $\langle  x_1^*,  x_2 \rangle = 0$. It follows from the generalized Perron theorem that $x_1^* \in {\rm int}(K_1^*)$. This inclusion implies that $(K_1 \cup (-
K_1)) \cap X_2' = \{0\}$. So we have that $x_2 \in {\rm
int}(T(K_1,K_2))\setminus (K_1 \cup (- K_1))$.

Let the statement of the theorem hold for $j-1$. Prove it for $j$.

Since $\rho(\wedge^j A)= \lambda_{1}\ldots \lambda_{j}$, the eigenvector
$\varphi_j$ corresponding to $\rho(\wedge^j A)$ can be represented
in the form of the exterior product $\varphi_j = x_1 \wedge \ldots \wedge x_j$ of the first $j-1$ eigenvectors $x_1 \in {\rm int}(K_1)$, $x_2 \in \widetilde{T}(K_1,K_2)$, $\ldots$, $x_{j-1} \in \widetilde{T}(K_1, \ \ldots, \
K_{j-1})$ and the $j$th
eigenvector $x_j$. The inclusion $\varphi_j \in {\rm int}(K_j)$ implies
$x_j \in \widetilde{T}(K_1, \ \ldots, \ K_{j}) \subseteq {\rm int}(T(K_1, \ \ldots, \ K_{j}))$.

 Let us show that the vector
$x_j$ does not belong to $T(K_1, \ \ldots, \
K_{j-1})$. Assume the contrary: let $x_j \in T(K_1, \ \ldots, \
K_{j-1})$. Since all the eigenvalues of the operator $A$ are distinct, it is not difficult to see, that $x_j$ belongs to the subspace $$X_j' = \{ x \in
{\mathbb R}^n: (c_1x_1^* + \ldots + c_{j-1}x_{j - 1}^*)(x) = 0 \quad \mbox{for any} \ c_1, \ldots, c_{j-1} \in {\mathbb R}\}.$$ Here $x_1^*, \ \ldots, \ x_{j-1}^*$ are the first $j-1$
eigenfunctionals of the operator $A$ corresponding to the
eigenvalues $\lambda_1, \ \ldots, \ \lambda_{j - 1}$, respectively. I.e. the vector $x_j$ is orthogonal to the subspace ${\rm Lin}(x_1^*, \ldots, x_{j - 1}^*)$ spanned by the vectors $x_1^*, \ldots, x_{j - 1}^*$. Since $x_j \in T(K_1, \ \ldots, \
K_{j-1}) \subseteq \widehat{T}(K_{j-1})$, we can find vectors $y_1, \ \ldots, \ y_{j-2}$ such that $y_1 \wedge \ldots \wedge y_{j-2} \wedge x_j \in (K_{j-1} \cup (- K_{j-1}))\setminus \{0\}$. Let us examine the exterior product $x_1^*\wedge \ldots \wedge x_{j - 1}^*$. It is not difficult to see, that $x_1^*\wedge \ldots \wedge x_{j - 1}^*$ belongs to ${\rm int}(K_{j-1}^*) \cup (- {\rm int}(K_{j-1}^*))$.
Examine the scalar product
$$\langle y_1 \wedge \ldots \wedge y_{j-2} \wedge x_j,  x_1^*\wedge \ldots \wedge x_{j - 1}^*\rangle = (y_1 \wedge \ldots \wedge y_{j-2} \wedge x_j)(x_1^*, \ldots, x_{j - 1}^*) = $$ $$ =\sum_{(i_{1}, \ldots, i_{j-1})} \chi(i_1, \ldots, i_{j-1}) \ \langle
y_1,x^*_{i_1} \rangle \ldots \langle y_{j-2},x_{i_{j-2}} \rangle \langle x_j,x^*_{i_{j-1}}\rangle = 0.$$
We came to the contradiction.

Finally, let us prove the inclusion $$\sum_{i=q}^p c_i x_i \in {\rm int}(T(K_1, \ \ldots, \
K_p))\setminus T(K_1, \ \ldots, \ K_{q - 1})$$ for each $1 \leq q \leq p
\leq n$, $c_p \neq 0$.

First let us prove, that any linear combination $\sum\limits_{i=q}^p c_i x_i$ belongs to \linebreak ${\rm int}(T(K_1, \ \ldots, \
K_p))$. Since the exterior product of the eigenvectors $x_1 \in {\rm int}(K_1)$, $x_2 \in \widetilde{T}(K_1,K_2)$, $\ldots$, $x_{p-1} \in \widetilde{T}(K_1, \ \ldots, \
K_{p-1})$, $x_p$ belongs to $ {\rm int}(K_p)$, we have $x_p \in \widetilde{T}(K_1, \ \ldots, \
K_p)$. Examine the exterior product
$x_1\wedge \ldots \wedge x_{p-1}\wedge (\sum\limits_{i=q}^p c_i x_i)$ which is obviously equal to $c_p(x_1\wedge \ldots \wedge x_{p-1}\wedge x_p)$. Since $c_p \neq 0$, we have $c_p(x_1\wedge \ldots \wedge x_{p-1}\wedge x_p) \in ({\rm int}(K_j)\cup {\rm int}(- K_j))$. So the inclusion $\sum\limits_{i=q}^p c_i x_i \in \widetilde{T}(K_1, \ \ldots, \
K_p) \subseteq {\rm int}(T(K_1, \ \ldots, \
K_p))$ is correct.

In the case when $c_p = 0$ and $\sum\limits_{i=q}^{p-1} c_i^2 \neq 0$, we construct a converging sequence $\{\sum\limits_{i=q}^{p-1} c_i x_i + \dfrac{1}{n}x_p\}_{n=1}^\infty \in {\rm int}(T(K_1, \ \ldots, \ K_p))$.

The fact, that $\sum\limits_{i=q}^p c_i x_i$ does not belong to $T(K_1, \ \ldots, \ K_{q -1})$ follows from the inclusion $\sum\limits_{i=q}^p c_i x_i \in X_q'$. As it is shown above, $X_q'$ has zero intersection with the set $T(K_{q - 1})$.
\end{pf}

Now we can state some necessary conditions for strict total positivity of a totally positive structure $\{K_1, \ \ldots, \ K_n\}$

\begin{cor}
Let a totally positive structure $\{K_1, \ \ldots, \ K_n\}$ on the space ${\mathbb R}^n$ be strictly totally positive. Then the intersection ${\rm int}(K_j) \cap \barwedge^j {\mathbb R}^n$ is nonzero for every $j \ (j = 2, \ \ldots, \ n)$.
\end{cor}
\begin{pf}
 If the structure $\{K_1, \ \ldots, \ K_n\}$ is strictly totally positive, then there exists at least one operator $A$ which is GSTP with respect to the structure
$\{K_1, \ \ldots, \ K_n\}$. Theorem 22 implies that the eigenvector $\varphi_j$ corresponding to the maximal eigenvalue of $\wedge^j A$, belongs to ${\rm int}(K_j)$. On the other hand, $\varphi_j$ is a simple $j$-vector, since it can be represented in the form of the exterior product $x_1\wedge \ldots \wedge x_j$, where $x_1, \ \ldots, \ x_j$ are the eigenvectors of the operator $A$.
\end{pf}

Note, that we can find a proper cone $K_j$ which belongs to $\wedge^j {\mathbb R}^n \setminus \barwedge {\mathbb R}^n$ for every $j = 2, \ \ldots, \ n-2$. Indeed, the set $\barwedge {\mathbb R}^n$ is closed in the space $\wedge^j {\mathbb R}^n$, so its complement $\wedge^j {\mathbb R}^n \setminus \barwedge^j {\mathbb R}^n$ is open.
Since $\wedge^j {\mathbb R}^n \setminus \barwedge {\mathbb R}^n$ is nonempty, we can find an element $\varphi \in (\wedge^j {\mathbb R}^n \setminus \barwedge {\mathbb R}^n)$. Since it is open, there exists a value $r > 0$ such that $B(\varphi, r) \subset (\wedge^j {\mathbb R}^n \setminus \barwedge {\mathbb R}^n)$ and the closed ball $\overline{B(\varphi, \dfrac{r}{2})} \subset B(\varphi, r) \subset (\wedge^j {\mathbb R}^n \setminus \barwedge {\mathbb R}^n)$. Thus we can construct a proper cone $K_j(\overline{B(\varphi, \dfrac{r}{2})})$ as it is shown in Example 3 (see Section 3).

\begin{cor}
Let the totally positive structure $\{K_1, \ \ldots, \ K_n\}$ on ${\mathbb R}^n$ be strictly totally positive. Then the set ${\rm int}(T(K_1, \ldots, K_j))$ is nonempty for every $j \ (j = 2, \ \ldots, \ n)$.
\end{cor}
\begin{pf}
 The proof obviously follows from Theorem 22.
\end{pf}

\begin{rmk}
It is not difficult to see, that the spectrum of a GTP operator is not neccesarily real, and not every GTP operator can be approximated by GSTP operators.
\end{rmk}

Analogically, it is not difficult to generalize Theorem 22 to the case of $k$-GSTP operators $(k = 2, \ \ldots, \ n)$ and GSSR operators.

Now we can state the following property of GSSR operators.

\begin{thm}
Let a linear operator $A: {\mathbb{R}}^{n} \rightarrow
{\mathbb{R}}^{n}$ be GSSR with respect to a totally positive structure $\{K_1, \ \ldots, \ K_n\}$. Then the interior of the set $T(K_j)$ is nonempty and the following inclusions hold for every $j = 1, \ldots, \ n$: $$A(T(K_j)) \subseteq {\rm int}(T(K_j)).$$
\end{thm}

\begin{pf}
It follows from the proof of Theorem 22 that a
GSSR operator $A$ has $n$ nonzero real simple eigenvalues $\lambda_1, \ \ldots, \ \lambda_n$ with the corresponding $n$ linearly independent eigenvectors $x_1, \ \ldots, \ x_n$. Moreover, we have the inclusion $x_j \in {\rm int}(T(K_j))$ for the $j$-th eigenvector $x_j$. Thus the interior of the set $T(K_j)$ contains at least one element for every $j = 1, \ldots, \ n$.

Let $x_1$ be an arbitrary vector from $T(K_j)$. Let us prove that
$Ax_1 \in {\rm int}(T(K_j))$. Since $T(K_j) \subseteq \widehat{T}(K_j)$, we can find the elements $x_2, \ \ldots, \ x_j \in {\mathbb{R}}^{n}$ such that $x_1 \wedge x_2 \wedge \ldots \wedge x_j \in (K_j\cup (-K_j)) \setminus \{0\}$. Examine the elements $Ax_2, \ \ldots, \ Ax_j$. Since the
operator $\wedge^j A$ is positive with respect to the cone $K_j$, we have the inclusion $Ax_1 \wedge Ax_2 \wedge \ldots \wedge Ax_j = (\wedge^j A)(x_1 \wedge x_2 \wedge \ldots \wedge x_j) \in ({\rm int}(K_j) \cup {\rm int}(-K_j)) \setminus \{0\}$. So $Ax_1 \in \widetilde{T}(K_j) = {\rm int}(T(K_j))$.
\end{pf}

\section{Classical theory of total positivity}
Let us examine the space ${\mathbb R}^n$ with the standard basis $e_1, \ \ldots, \ e_n$ and its $j$th $(j = 2, \ \ldots, \ n)$ exterior powers $\wedge^j {\mathbb R}^n$ with the canonical basis which consists of the exterior products of the
form $\{e_{i_1}\wedge \ldots \wedge e_{i_j}\}$ where $1 \leq i_1
< \ldots < i_j \leq n$. As it is mentioned above, we denote the cone spanned by basic vectors by ${\mathbb R}^n_+$ and the cone spanned by the $j$th exterior basic vectors by $\wedge^j {\mathbb R}^n_+$. $S^-(x)$ denotes the number of sign changes in the sequence $(x^1, \ \ldots, \ x^n)$ of the coordinates with zero terms discarded. $S^+(x)$ denotes the maximum number of sign changes in the sequence $(x^1, \ \ldots, \ x^n)$ where zero terms are arbitrarily assigned values $\pm 1$.

The following lemma describes the link between the sign changes of vectors in ${\mathbb R}^n$ and their exterior products (see \cite{AN}, p. 198, Lemma 5.1).

\begin{lem}
 Let $x_1, \ldots, x_j \in {\mathbb R}^n$. In order for $$ S^+(\sum_{i=1}^j c_i x_i) \leq j-1$$
for each $c_1, \ \ldots, c_j \in {\mathbb R}$ and $\sum\limits_{i=q}^p c_i^2 \neq 0$ it is neccesary and sufficient that $x_1\wedge \ldots \wedge x_j \in ({\rm int}(\wedge^j {\mathbb R}^n_+) \cup {\rm int}(- \wedge^j {\mathbb R}^n_+))$.
\end{lem}

Examine the set $$M(j) = \{x \in {\mathbb R}^n; \ S^-(x) \leq j-1\}.$$  The set $M(j)$ is closed, solid and uniform (see, for example, \cite{KRASN1}, \cite{SOB1}). The following equality for its interior is valid: $${\rm int}(M(j)) = \{x \in {\mathbb R}^n; \ S^+(x) \leq j-1 \}.$$

\begin{prop}
The following equalities hold for every $j = 2, \ \ldots, \ n$:
 $$M(j) = T(\wedge^j {\mathbb R}^n_+) = T({\mathbb R}^n_+, \ldots, \wedge^j {\mathbb R}^n_+);$$
$${\rm int}(M_j) = \widetilde{T}(\wedge^j {\mathbb R}^n_+) = \widetilde{T}({\mathbb R}^n_+, \ldots, \wedge^j {\mathbb R}^n_+).$$
\end{prop}
\begin{pf}
The proof obviously follows from Lemma 26.
\end{pf}

Let us recall the following definitions. A matrix ${\mathbf A}$ of a linear operator $A:{\mathbb R}^n \rightarrow {\mathbb R}^n$ is called {\it totally positive (TP)} if it is nonnegative and its $j$-th compound matrices ${\mathbf A}^{(j)}$ are also nonnegative for every $j$ $(j = 2, \ \ldots, \ n)$.

A matrix ${\mathbf A}$ of a linear operator $A:{\mathbb R}^n \rightarrow {\mathbb R}^n$ is called {\it strictly totally positive (STP)} if it is positive and its $j$-th compound matrices ${\mathbf A}^{(j)}$ are also positive for every $j$ $(j = 2, \ \ldots, \ n)$.

A matrix ${\mathbf A}$ of a linear operator $A:{\mathbb R}^n \rightarrow {\mathbb R}^n$ is called {\it sign-regular (SR)}, if there exist numbers $\epsilon_1, \ \ldots, \ \epsilon_n$ each equal to $\pm 1$ such that the matrices
$\epsilon_j{\mathbf A}^{(j)}$ are nonnegative for every $j$ $(j = 1, \ 2, \ \ldots, \ n)$.

A matrix ${\mathbf A}$ of a linear operator $A:{\mathbb R}^n \rightarrow {\mathbb R}^n$ is called {\it strictly sign-regular (SSR)} if there exist such numbers $\epsilon_1, \ \ldots, \ \epsilon_n$ each equal to $\pm 1$ such that the matrices
$\epsilon_j{\mathbf A}^{(j)}$ are positive for every $j$ $(j = 1, \ 2, \ \ldots, \ n)$.

Now Theorem 25 turns into the following statement (see, for example, \cite{PINK}).

\begin{thm}
 Let the matrix ${\mathbf A}$ of a linear operator $A:{\mathbb R}^n \rightarrow {\mathbb R}^n$ be SSR. Then the following inequality holds for each nonzero vector $x \in {\mathbb R}^n$:
$$S^+(Ax) \leq S^-(x).$$
\end{thm}

Theorem 22 turns into the classical Gantmacher--Krein theorem (Theorem 2).

\section{Totally J-sign-symmetric matrices}

Since it is easy to see if a matrix is diagonally similar to a nonnegative one, let us reformulate the given above definitions and theorems in terms of compound matrices. In this case the conditions of generalized total positivity become easily verified.

Let us give the following definitions.

A matrix ${\mathbf A}$ of a linear operator $A:{\mathbb R}^n \rightarrow {\mathbb R}^n$ is called {\it totally J--sign-symmetric (TJS)}, if it is J--sign-symmetric, and its $j$-th compound matrices ${\mathbf A}^{(j)}$ are also J--sign-symmetric for every $j$ $(j = 2, \ \ldots, \ n)$.

A matrix ${\mathbf A}$ of a linear operator $A:{\mathbb R}^n \rightarrow {\mathbb R}^n$ is called {\it strictly totally J--sign-symmetric (STJS)}, if it is strictly J--sign-symmetric, and its $j$-th compound matrices ${\mathbf A}^{(j)}$ are also strictly J--sign-symmetric for every $j$ $(j = 2, \ \ldots, \ n)$.

It is easy to see, that the class of TP matrices belongs to the class of TJS matrices, and the class of STP matrices belongs to the class of STJS matrices.

Note that it is also not difficult to reformulate the definitions of generalized sign-regularity and generalized strict sign-regularity in terms of compound matrices and to introduce the classes of $k$-totally J--sign-symmetric and strictly $k$-totally J--sign-symmetric matrices for every $k = 2, \ \ldots, \ n$.

Now we examine basic properties of TJS and STJS matrices.

\begin{prop}
Let the matrix ${\mathbf A}$ of a linear operator $A: {\mathbb R}^n \rightarrow {\mathbb R}^n$ be TJS (STJS). Then ${\mathbf A}^T$ (the transpose of ${\mathbf A}$), as well as every principal submatrix of ${\mathbf A}$ and ${\mathbf A}^T$ is TJS (respectively STJS).
\end{prop}

\begin{pf}
As $M$ denotes an arbitrary subset of $[n]$, ${\mathbf A}(M)$ is the principal submatrix which consists of the rows and columns with the numbers from $M$. Let us consider ${\mathbf A}(M)$ as the matrix of the restriction $A|_{L(M)}$ of $A$ to the subspace $L(M)$ spanned by the basic vectors with the numbers from $M$. Since ${\mathbf A}$ is JS, the operator $A$ leaves invariant some basic cone $K \subset {\mathbb R}^n$. It is not difficult to see, that $A|_{L(M)}$ leaves invariant the set $K(M) = {\rm pr}_{L(M)}K$ i.e the projection of $K$ on the subspace $L(M)$. According to Property 1 of basic cones (see Section 3), the set $K(M)$ is a basic cone in the space $L(M)$. So we conclude that the submatrix ${\mathbf A}(M)$ is JS.

Applying the same reasoning to the $l$th compound matrix
$({\mathbf A}(M))^{(l)}$ (here $l = 2, \ \ldots, \ {\rm Card}(M)$), we obtain that $({\mathbf A}(M))^{(l)}$ is also JS.
The case of STJS matrices is considered analogically.
The fact, that ${\mathbf A}^T$, as well as every principal submatrix of ${\mathbf A}^T$ is TJS (STJS) follows from Proposition 14.
\end{pf}
\begin{rmk}
 If ${\mathbf A}$ is a TP (STP) matrix then every submatrix of ${\mathbf A}$ and ${\mathbf A}^T$ is obviously TP (respectively, STP). However, if ${\mathbf A}$ is a TJS (STJS) matrix then the analogous of this statement is true only for principal submatrices. It is easy to see, that an arbitrary submatrix of a TJS  matrix ${\mathbf A}$ may not be JS.
\end{rmk}

Since the projection of an arbitrary proper cone on a basic subspace may not be a proper cone in this subspace, Proposition 29 may not be correct for an arbitrary GTP (GSTP) operator.

{\bf Example.}
Let a linear operator $A: {\mathbb R}^3 \rightarrow {\mathbb R}^3$ be a rotation operator with the angle $\theta$ around the axis defined by $e_3$. It has the following matrix representation:

$${\mathbf A} = \begin{pmatrix}
  \cos(\theta) & -\sin(\theta) & 0 \\
\sin(\theta) & \cos(\theta) & 0 \\
 0 & 0 & 1
\end{pmatrix}.$$

It is not difficult to see, that the rotation operator $A$ leaves invariant the ice-cream cone $K_1$ defined in the following way: $$K_1 = \{x = (x^1, x^2, x^3) \in {\mathbb R}^3; \sqrt{(x^1)^2+ (x^2)^2}\} \leq x^3\}.$$

Examine the second exterior power $\wedge^2 A$ of the operator $A$. It is represented by the second compound matrix ${\mathbf A}^{(2)}$ in the basis $e_1
\wedge e_2$, $e_1 \wedge e_3$, $e_2 \wedge e_3$.

$${\mathbf A}^{(2)} = \begin{pmatrix}
1 & 0 & 0 \\
0 & \cos(\theta) & -\sin(\theta) \\
0 & \sin(\theta) & \cos(\theta)
\end{pmatrix}.$$

It is obvious, that $\wedge^2 A$ is a rotation in $\wedge^2 {\mathbb R}^3$ with the same angle $\theta$ around the axis defined by $e_1 \wedge e_2$. So it leaves invariant the ice-cream cone $K_2$ defined in the following way: $$K_2 = \{\varphi = (\varphi^1, \varphi^2, \varphi^3) \in \wedge^2 {\mathbb R}^3; \sqrt{(\varphi^2)^2+ (\varphi^3)^2}\} \leq \varphi^1\}.$$

The matrix of the third exterior power $\wedge^3 A$ consists of only one element $\det{\mathbf A} = 1$. So $\wedge^3 A$ leaves invariant the cone $K_3$ which is the positive real axis.

Thus the operator $A$ is GTP with respect to the structure $\{K_1, \ K_2, \ K_3\}$.

But the principal submatrix ${\mathbf A}(1,2)$ which is the matrix of the restriction of $A$ on a basic subspace $L(1,2)$ spanned by the vectors $e_1, e_2$ is not GTP. It is not $K$-nonnegative. Indeed, examine the eigenvalues of
$${\mathbf A}(1,2) = \begin{pmatrix}
  \cos(\theta) & -\sin(\theta)  \\
\sin(\theta) & \cos(\theta)
\end{pmatrix}.$$
In the case of $0< \theta < \pi$ they are $\lambda_1 = \cos(\theta) + i \sin(\theta)$, $\lambda_2 = \cos(\theta) - i \sin(\theta)$ and both are complex.

The reasoning of the proof of Proposition 29 is not valid in this case, since the projection of the ice-cream cone $K_1$ on $L(1,2)$ coincides with the whole $L(1,2)$.

\begin{prop}
 Let the matrix ${\mathbf A}$ of a linear operator $A: {\mathbb R}^n \rightarrow {\mathbb R}^n$ be TJS (STJS). Let $P$ be an arbitrary permutation of $[n]$ and ${\mathbf P}$ be the corresponding permutation matrix. Then the matrix ${\mathbf P}{\mathbf A}{\mathbf P}^{-1}$ obtained from the initial matrix by re-numerating of both the rows and columns in order $P$ is TJS (respectively, STJS). In particular, the matrix obtained from $A$ by reversing the order of both its rows and columns is TJS (STJS) of the same structure.
\end{prop}

\begin{pf}
The proof is obvious since any similarity transformation with the permutation matrix is just re-numbering of the basic vectors.
\end{pf}

\begin{thm} Let the matrix ${\mathbf A}$ of a linear
operator $A: {\mathbb R}^n \rightarrow {\mathbb R}^n$ be
STJS. Then
all the eigenvalues of the operator $A$ are positive and simple:
$$\rho(A) = \lambda_1 > \lambda_2 > \ldots > \lambda_n > 0.$$ \end{thm}

\begin{cor}
 All principal minors of a STJS matrix ${\mathbf A}$ are positive.
\end{cor}
\begin{pf}
The positivity of every real eigenvalue of ${\mathbf A}$ as well as of every real eigenvalue of each principal submatrix of ${\mathbf A}$ implies the positivity of all principal minors of ${\mathbf A}$ (see \cite{FIP}, p. 385, Theorem 3.3).
\end{pf}

\section{Conclusions}
Many important properties of GTP and GSTP operators like the criteria of generalized total positivity and generalized strict total positivity, the factorization of GTP and GSTP operators, determinantal inequalities, the interlacing properties of the eigenvalues as well as many important examples of GTP operators are not considered in this paper. The application of the obtained theory to multi-dimensional boundary-value problems is also not mentioned. However, the author hopes that it would be possible to state GSTP properties of the corresponding Green's functions for certain classes of such problems. This would imply the positivity of the spectra of the corresponding differential operators.



\bibliographystyle{elsarticle-num}
\bibliography{<your-bib-database>}







\end{document}